\documentclass[11pt]{amsart}
\usepackage{amssymb}

\numberwithin{equation}{section}

\newcommand{\req}[1]{(\ref{#1})}

\def\b{{\beta}}
\def\Q{{\Bbb Q}}

\def\N{{\Bbb N}}
\def\bJ{{\bf J}}
\def\R{{\Bbb R}}
\def\P{{\Bbb P}}
\def\E{{\Bbb E}}

\def\BB{{\bf B}}

\def\BB{{\bf B}}

\def\BJ{{\bf J}}

\def\BR{{\bf R}}
\def\BD{{\bf D}}
\def\BJ{{\bf J}}

\def\bmu{{\boldsymbol \mu}}
\def\wka{{\widetilde \kappa}}
\def\bq{{\bf q}}
\def\bb{{\bf b}}
\def\ba{{\bf a}}
\def\Bx{{\bf x}}
\def\bx{{\bf x}}
\def\buu{{\bf u}}
\def\bv{{\bf v}}
\def\wbx{\tilde{\bf x}}
\def\by{{\bf y}}
\def\By{{\bf y}}
\newcommand{\reals}{{\Bbb{R}}}

\def\bx{\Bx}

\def\CL{{\mathcal L}}

\def\CL{{\mathcal L}}

\def\Aa{{\mathcal A}}

\def\Ca{{\mathcal C}}

\def\Ea{{\mathcal E}}
\def\Fa{{\mathcal F}}
\def\Gga{{\mathcal G}}
\def\Ha{{\mathcal H}}

\def\C{{\Bbb C}}
\def\R{{\Bbb R}}
\def\Q{{\Bbb Q}}
\def\N{{\Bbb N}}
\def\E{{\Bbb E}}

\def\half{\frac{1}{2}}
\def\oneN{\frac{1}{N}}
\def\prf{\noindent {\bf Proof:\ }}

\def\nn{\noindent}

\def\pm{{probability measure\ }}

\def\a{\alpha}
\def\b{\beta}
\def\d{\delta}
\def\e{\epsilon}
\def\g{\gamma}

\def\s{\sigma}
\def\t{\theta}

\def\NC{{\mbox{NC}}}
\def\D{\Delta}

\def\L{\Lambda}

\def\part{\partial}

\def\ts{\times}

\def\ra{\rightarrow}

\def\lbc{\lbrace}
\def\rbc{\rbrace}

\def\tilde{\widetilde}

\cleardoublepage
\newtheorem{prop}{Proposition}[section]

\newtheorem{lem}[prop]{Lemma}

\newtheorem{cor}[prop]{Corollary}
\newtheorem{theo}[prop]{Theorem}
\newtheorem{ass}[prop]{Hypothesis}

\title[Dynamics for spherical p-spins]
{Cugliandolo-Kurchan equations for dynamics of Spin-Glasses.}
\author{Gerard Ben Arous}
\address{Courant Institute\\
251 Mercer Street\\
New York, NY 10012\\
and Mathematics, EPFL\\
1015 Lausanne, Switzerland}

\email{benarous@cims.nyu.edu}

\author{Amir Dembo}
\address{Department of Statistics and Department of Mathematics\\
Stanford University\\ Stanford, CA 94305.}
\email{amir@math.stanford.edu}

\author{Alice Guionnet}

\thanks{\noindent Research partially supported by NSF
grants \#DMS-0406042, \#DMS-FRG-0244323.
\newline
{\bf AMS (2000) Subject Classification:}
{Primary: 82C44 Secondary:  82C31, 60H10, 60F15, 60K35}
\newline
{\bf Keywords:} Interacting random processes, Disordered systems,
Statistical mechanics, Langevin dynamics, Aging, $p$-spin models.}
\address{UMPA, Ecole Normale Superieure de Lyon\\
46 all\'ee d'Italie\\
69364 Lyon Cedex 07, France}
              
\email{Alice.Guionnet@umpa.ens-lyon.fr}

\date{September 3, 2004}
              
\begin{document}

\begin{abstract}
We study the Langevin dynamics for the family of
spherical $p$-spin disordered mean-field models 
of statistical physics. We prove that in the 
limit of system size $N$ approaching infinity, 
the empirical state correlation and
integrated response functions for these 
$N$-dimensional coupled diffusions converge
almost surely and uniformly in time, 
to the non-random 
unique strong solution of a pair of explicit non-linear 
integro-differential equations,
first introduced by Cugliandolo and Kurchan.    
\end{abstract}

\maketitle

{\vspace{1cm}}
\nn
{\underline{\bf Keywords :}}
Interacting random processes, Disordered systems,
 Statistical mechanics,
Langevin dynamics.

\nn
{\underline{\bf  Mathematics Subject of Classification  :}}
60H10, 82B44, 60K35,  82C31, 82C22.

\section{Introduction and main results} 

Markovian dynamics with random interactions can produce very complex
  phase transitions, and fascinating long time behaviors, for 
strong disorder (or low temperature). The physics 
literature has shown that dynamics of mean-field spin glasses 
is a very good field to get an accurate sample of possible and 
generic long time phenomena (as aging, memory, rejuvenation, 
failure of the Fluctuation-Dissipation theorem, 
see \cite{BKM} for a good survey). 

\smallskip
This class of problems can be roughly described as follows.
Let $\Gamma$ (a compact metric space) be the state space 
for spins and $\nu$ be a probability measure on $\Gamma$.
Typically, in the discrete or Ising spins context $\Gamma=\{-1,1\}$ 
and $\nu={1/2}(\delta_{1}+\delta_{-1})$. 
In the continuous or soft spin context, $\Gamma=I$ 
a compact interval  of the real line and 
$\nu(dx) = Z^{-1}e^{-U(x)}dx$ , where $U(x)$ is the "one-body potential". 
For each configuration of the spin system, i.e. 
for each  $\bx=(x_1,...,x_N) \in \Gamma^N$ one defines a random Hamiltonian, 
$H^N_{\BJ}(\bx)$, as a function of the configuration $\bx$ 
and of an exterior source of randomness $\BJ$, i.e. a random variable 
defined on another probability space.
The Gibbs measure at inverse temperature $\b$ is 
then defined on the configuration space $\Gamma^N$ by
$$\mu^N_{\b,\BJ}(d\bx)= \exp({-\b H^N_{\BJ}(\Bx))\nu(d\bx)}/Z^N_{\BJ}$$
The statics problem amounts to understanding the large $N$
behavior of these measures for various classes of random Hamiltonians (\cite{Talagrand} is a recent and beautiful book on the mathematical results pertaining to these equilibrium problems).
The dynamics question consists of understanding the behavior of 
Markovian processes on the configuration space $\Gamma^N$, for which 
the Gibbs measure is invariant and even reversible, 
in the limit of large systems (large $N$) and long times, 
either when the randomness $\BJ$ is fixed (the quenched case) 
or when it is averaged (often called the annealed case, 
in the mathematics literature, but not in the physics papers).
 
These dynamics are typically Glauber dynamics 
for the discrete spin setting, or Langevin dynamics for continuous spins. 
Defining precisely what we mean here by 
large system size and long time is a very important question, 
and very different results can be expected (and sometimes proved, see \cite{ICM} and references therein),
for various time-scales as functions of the size of the system. 
We will restrict ourselves to the case where the system size
is first taken to infinity, i.e for the shortest possible long time scales, 
much too short typically to allow any escape from meta-stable states (as opposed to the situation in \cite{BBG} for instance).
The first step is then to derive limiting 
equations for various quantities, when $N$ tends to $\infty$. 
The second step is to understand the large time behavior of 
these limiting macroscopic equations.
Let us be more specific by describing one of the main initial 
objectives of the theory, which is 
the long time behavior of the Langevin dynamics 
of the Sherrington-Kirkpatrick model and its generalization, the $p$-spin model. 
In the SK model, either for discrete or for continuous spins, the 
Hamiltonian is given by:
$$H^N_{\BJ}(\Bx)= \sum_{1\le i,j \le N} J_{\{ij\}} x^{i}x^{j}\,,
$$
where the randomness is due to the couplings $J_{\{ij\}}$ 
which are assumed to be i.i.d Gaussian centered, of variance $N^{-1}$
in case $i \neq j$ and $2N^{-1}$ in case $i=j$.
In the $p$-spins model the Hamiltonian contains 
interaction between subsets of spins of size $p$, that is,
\begin{equation}\label{pspins}
H^N_{\BJ}(\Bx)= \sum_{1\le i_1 < \ldots< i_p \le N} 
J_{i_1\ldots i_p} x^{i_1}\ldots  x^{i_p}\,,
\end{equation}
where the couplings $J_{i_1\ldots i_p}$ are assumed to be i.i.d. 
Gaussian centered and of variance of $O(N^{-(p-1)})$.
The initial SK model corresponds of course to the case $p=2$.

Propagation of chaos for dynamics of the SK model, 
or equivalently the large $N$ limit for the behavior 
of the empirical measure $\frac{1}{N} \sum_{i=1}^N \delta_{x^i(t)}$ 
has been obtained some time ago both for continuous spins (see 
\cite{MPV,SZ},
in the physics literature and \cite{BAG1,BAG2,G1} for a mathematical treatment) 
or discrete ones (see \cite{MG1}).
The limiting equations (the so called self-consistent single-spin dynamics) 
are very complex and have resisted so far all 
attempts to understand their long time behavior. This is due in part 
to the fact that the empirical measure is a much too rich object. 
Finding an autonomous system of tractable equations for a 
proper well chosen set of lower dimensional quantities is 
still a very open question, even in the physics literature.

But a large range of interesting and related models have been 
recently analyzed more successfully in the physics literature, 
i.e the spherical $p$-spin models (see \cite{BKM,LesHouches}). 
The spherical version of the spin glass models consists of  
a classical simplification, that is, 
replacing the product structure 
of the configuration space $\Gamma^N$ by the sphere $S^{N-1}({\sqrt{N}})$ 
of radius  $\sqrt{N}$ in $\R^N$, in effect,  
imposing on the configuration $\bx$ 
a hard constraint $\frac{1}{N}\sum_{i=1}^N x_i^2 = 1$. 
The spherical $p$-spin Gibbs measure is thus the 
probability measure on the sphere $S^{N-1}({\sqrt{N}})$ given by 
$$\mu^N_{\b,\BJ}(d\bx)= \exp( -\b H^N_{\BJ}(\Bx))\nu_N (d\bx) /Z^N_{\BJ}\,,
$$
where the Hamiltonian is given by \req{pspins} and the measure $\nu_N$ 
is the uniform measure on the sphere $S^{N-1}({\sqrt{N}})$.
One can also study a very similar problem by replacing the hard 
spherical constraint by a soft one, i.e by replacing the uniform 
measure $\nu_N$  on the sphere $S^{N-1}({\sqrt{N}})$ by a measure on $\R^N$
$$
\nu_N(d\bx) = \exp{\big(- N f(\frac{1}{N}\sum_{i=1}^N x_i^2)\big)} 
d\bx /Z^N\,,
$$
where $f$ is a smooth function growing fast enough at infinity.
This study has been done successfully (see \cite{BDG}, and \cite{CD}) 
in the case where $p=2$, the spherical SK model.
There one could obtain a very complete description 
of the limiting dynamics using only one quantity, 
the empirical state correlation function
$$ C_N(s,t)= \frac{1}{N}\sum_{i=1}^N x^i(s)x^i(t) \,. $$
Indeed, this quantity was shown to have a non-random
limit $ C(s,t)$, when $N$
tends to $\infty$, which satisfies an autonomous integro-differential 
equation.
This is a rather simple setting, when $p=2$, since this quadratic 
case can use the well know tools pertaining to the study 
of spectra of GOE random matrices, after diagonalization 
of the random matrix $ (J_{ij})$. The results, though showing
an aging phenomenon, are very far from what is expected for 
the true (i.e non spherical) SK model.
The case of $p>2$ is a completely different story. 
The physics literature, mainly Cugliandolo and Kurchan 
(see \cite{BKM,LesHouches,CK}), 
has given a very rich picture of the behavior of the Langevin dynamics 
for these spherical $p$-spin models, which are believed to mimic in certain 
cases the behavior of the dynamics of the true SK model. The main innovation 
is the fact that the description of the limiting dynamics relies now on 
two objects, the empirical correlation function and 
the response function (and not only one as in the spherical $p=2$
 case or the whole infinite dimensional object of the true, non-spherical,
 SK model). These two functions are believed to satisfy a 
set of coupled non-linear integro-differential equations 
(which of course decouple in the $p=2$ case).
Our work establishes rigorously the asymptotic validity of these 
Cugliandolo-Kurchan equations 
(see \cite{CK} for the physical derivation of these).
We shall devote a future work to the fascinating question 
of understanding the large time behavior of the solution to these equations
(see \cite{BKM} or \cite{LesHouches} to see what is expected).  
Some of the needed mathematical tools for 
attacking this problem 
can be found in \cite{GM}, as well as an
intriguing link with equations arising naturally in non-commutative 
probability.

We now describe more precisely our result and then the structure of this
paper.

Fixing a positive integer $N$ (denoting the system size), we
consider the mean field random Hamiltonian
\begin{equation}\label{eq:hamil}
H^N_{\BJ}(\Bx)=-2 \sum_{p=1}^m \frac{a_p}{p!}
 \sum_{1\le i_1,\ldots,
i_p \le N} J_{i_1\ldots i_p} x^{i_1}\ldots  x^{i_p} \,,
\end{equation}
where $m \geq 2$, the state variable is
$\Bx=(x^{1}, \ldots ,x^{N}) \in \reals^N$,
and the disorder parameters $J_{i_1 \ldots i_p}=J_{\{i_1,\ldots,i_p\}}$
are independent (modulo the permutation of the indices) centered
Gaussian variables. The variance of $J_{i_1\ldots i_p}$ is
$c(\{i_1,\ldots,i_p\}) N^{-p+1}$, where
\begin{equation}\label{eq:vardef}
c(\{i_1,\ldots,i_p\})=\prod_k l_k! \,,
\end{equation}
and $(l_1,l_2,\ldots)$ are
the multiplicities of the different elements of the set      
$\{i_1,\ldots,i_p\}$ (for example, $c=1$ when $i_j \neq i_{j'}$ for
any $j \neq j'$, while $c=p!$ when all $i_j$ values are the same).

Let $f$ be a differentiable function on $\reals_+$ with
$f'$ locally Lipschitz, such that 
\begin{equation}\label{eq:fcondu}
\sup_{\rho \geq 0} |f'(\rho)| (1+\rho)^{-r} < \infty
\end{equation}
for some $r<\infty$, and for some $A,\d>0$,
\begin{equation}\label{eq:fcond}
\inf_{\rho \geq 0} \{ f'(\rho) - A \rho^{m/2+\d-1} \} > -\infty
\end{equation}
(typically, $f(\rho)=\kappa (\rho-1)^r$ for some $r>m/2$
and $\kappa \gg 1$).
At temperature $\b^{-1}>0$, the soft spherical version of the system
corresponds to the equilibrium probability measure
$\mu^N_{\b,\BJ}$ on $\reals^N$ whose density
(with respect to Lebesgue measure) is
$$
\frac{d \mu^N_{\b,\BJ}}{d\Bx} = Z_{\b,\BJ}^{-1}
e^{-\b H_{\BJ}^N (\Bx) - N f(N^{-1} \|\Bx\|^2)} \,,
$$
where $\|\cdot\|$ denotes the Euclidean norm on $\R^N$ and
the normalization factor
$Z_{\b,\BJ}= \int
e^{-\b (H_{\BJ}^N (\Bx) -N f(N^{-1} \|\Bx\|^2)} d\Bx$ is a.s.
finite
(by (\ref{eq:fcond})).
Recall that $\mu^N_{\b,\BJ}$ is the invariant measure of the
randomly interacting particles described by the
(Langevin) stochastic differential system:
\begin{equation}\label{interaction}
dx^j_t=dB_t^j- f'(N^{-1} \|\Bx_t\|^2) x^j_t dt
+\beta G^j(\Bx_t) dt \,,
\end{equation}
where $\BB=(B^1,\ldots,B^N)$ is an N-dimensional standard
Brownian motion, independent of both the initial condition $\Bx_0$ and
the disorder $\BJ$, while
\begin{equation}\label{eq:gdef}
G^i(\Bx):=-\frac{1}{2} \partial_{x^i}
\Big( H^N_{\BJ} (\Bx) \Big) =\sum_{ p=1}^m \frac{a_p}{(p-1)!}
 \sum_{1\le i_1,\ldots, i_{p-1} 
\le N} J_{i i_1\ldots i_{p-1}} x^{i_1}\ldots  x^{i_{p-1}} \,,
\end{equation}
for $i=1,\ldots,N$.

In Proposition \ref{existN}
we prove that for a.e. disorder
$\BJ$, initial condition  $\Bx_0$ and Brownian path $\BB$,
there exists a unique strong solution of \req{interaction} for
all $t \geq 0$, whose law we denote by $\P^{N}_{\b,\Bx_0,{\BJ}}$.

We are interested in the time
evolution of the empirical state correlation function
\begin{equation}\label{eq:cdef}
C_N(s,t):=\oneN \sum_{i=1}^N x^i_s x^i_t\,,
\end{equation}
and that of the empirical integrated response function
\begin{equation}\label{eq:chidef}
\chi_N(s,t):=\oneN \sum_{i=1}^N x^i_s B^i_t\,,
\end{equation}
under the quenched  law $\P^{N}_{\b,\Bx_0,\BJ}$, as the
system size $N \to \infty$.
Note that since $\{\chi_N(s,t),0\le t\le s\le T\}$
is not determined by $N^{-1}\sum\d_{x_{[0,T]}^i}$, 
the strategy by which the
limiting equations in \cite{BAG1,BAG2} are derived,
does not apply directly in our case. 

We assume hereafter that the initial condition $\Bx_0$ is
independent of the disorder $\BJ$, and the limit
\begin{equation}\label{eq:x0cond}
\lim_{N \to \infty} \E C_N(0,0) = C(0,0) \,,
\end{equation}
exists, and is a finite. Further, we assume that 
the tail probabilities $\P(|C_N(0,0)-C(0,0)|>x)$ decay
exponentially fast in $N$ (so the convergence $C_N(0,0) \to C(0,0)$
holds almost surely), and that
for each $k<\infty$, the sequence $N \mapsto \E [ C_N(0,0)^k ]$ is
uniformly bounded.

To be more specific, we 
consider hereafter the product probability space $\Ea_N
 = \R^N  \times \R^{d(N,m)}\times \C([0,T],\R^N) $
(here $T$ is a fixed time and $d(N,m)$
is the dimension of the space of the interactions $\bJ$), 
equipped with the natural Euclidean norms for the finite dimensional parts, 
i.e $(\Bx_0, \BJ)$, and the sup-norm for the Brownian motion $\BB$. 
The space $\Ea_N$ is endowed with 
the product probability measure $\P = \mu_N \otimes \gamma_N \otimes P_N$, 
where  $\mu_N$ denotes the distribution of $\Bx_0$, $\gamma_N$ 
is the (Gaussian) distribution of the coupling constants $\BJ$, 
and $P_N$ is the distribution of the $N$-dimensional Brownian motion. 
\begin{ass}\label{defconcentration}
For  $(\Bx_0, \BJ, \BB) \in \Ea_N$ we introduce the norms
$$\|(\Bx_0, \BJ, \BB)\|^2= \sum_{i=1}^N(x_0^i )^2 
+ \sum_{p=1}^m \sum_{1 \leq i_1\ldots i_{p} \leq N} 
(N^{\frac{p-1}{2}}J_{i_1\cdots i_{p}})^2
+ \sup_{0 \leq t \leq T}\sum_{i=1}^N(B_t^i)^2\,.$$
We shall assume that $\mu_N$ is such that 
 the following concentration of measure property 
holds 
on $\Ea_N$;   there exists two finite positive 
constants $C$ and $\alpha$, independent
on $N$, such that, if $V$
is a Lipschitz function on $\Ea_N$, with Lipschitz constant $K$, then
for all $\rho>0$,
$$\mu_N \otimes \gamma_N \otimes P_N [ |V-\E[V]|\geq \rho ] 
\leq C^{-1} \exp{(-C(\frac{\rho}{K})^\alpha)}\,.$$
\end{ass}

The above concentration inequality holds for 
any Lipschitz function $V$ that does not depend on $\Bx_0$,
since the concentration of measure property holds for the 
Gaussian measures $\gamma_N \otimes P_N$, with $\alpha=2$
(c.f. \cite{lesToulousains}).  Unfortunately, assuming that 
the concentration of measure property holds for the measures $\mu_N$
on $\R^N$ does not assure the concentration 
of measure for the product measure $\mu_N \otimes \gamma_N \otimes P_N$. 
Hence, 
we shall have to assume a property which implies the 
concentration inequality of Hypothesis 
 \ref{defconcentration}
and further tensorizes 
(see \cite{lesToulousains} for a more thorough discussion).
For example, if $\mu_N$ satisfy the 
Poincar\'{e} inequality uniformly in $N$,
then the product measures $\mu_N \otimes \gamma_N \otimes P_N$ also 
satisfy the Poincar\'{e} inequality uniformly in $N$, 
since the Gaussian measure $\gamma_N \otimes P_N$ does. 
The required uniform in $N$ 
concentration property is then satisfied, with $\alpha=1$
(c.f. \cite{AidaStroock,lesToulousains}). 

%

Our main result is the proof that as $N \to \infty$
the functions $C_N(s,t)$ and $\chi_N(s,t)$ converge to
non-random functions $C(s,t)$ and $\chi(s,t)$, that are
characterized by the following theorem.
\begin{theo}\label{thm-macro}
Let $\psi(r)=\nu'(r)+r\nu''(r)$ and
\begin{equation}\label{eq:nudef}
\nu(r):=\sum_{p=1}^m \frac{a_p^2}{p!} r^p \,.
\end{equation}
Suppose $\mu_N$ satisfies hypothesis \ref{defconcentration}. 
Fixing any $T<\infty$, as $N \to\infty$ the
random functions $C_N$ and $\chi_N$ converge uniformly on $[0,T]^2$,
almost surely and in $L_p$ with respect to $\Bx_0$, $\BJ$ and $\BB$
to non-random functions $\chi(s,t)=\int_0^t R(s,u) du$ and $C(s,t)=C(t,s)$.
Further, $R(s,t)=0$ for $t>s$, $R(s,s)=1$,
and for $s>t$ the absolutely continuous functions 
$C$, $R$ and $K(s)=C(s,s)$
are the unique solution in the space of bounded, continuous functions,
of the integro-differential equations
\begin{eqnarray}
\partial_s R(s,t)&\!\!\!\! =\!\!\!\!
& - f'(K(s)) R(s,t) + \b^2 \int_t^s
R(u,t) R(s,u) \nu''(C(s,u)) du ,\label{eqR}\\
\partial_s C(s,t)&\!\!\!\! = \!\!\!\!
& - f'(K(s)) C(s,t) +
\b^2 \int_0^s C(u,t) R(s,u) \nu''(C(s,u)) du
+ \b^2 \int_0^t \nu'(C(s,u)) R(t,u) du,\label{eqC}\\
\partial_s K(s) &\!\!\!\! =\!\!\!\! & -2 f'(K(s)) K(s) + 1 + 2\b^2
\int_0^s \psi(C(s,u)) R(s,u) du , \label{eqZ}
\end{eqnarray}
where the initial condition
$K(0)=C(0,0)$ is determined by (\ref{eq:x0cond}).
\end{theo}

We note in passing that for $p=2$, i.e. $\nu(r)=r^2/2$, we get from \req{eqR}
the autonomous equation 
$$
\partial_s H(s,t) =\b^2 \int_t^s H(u,t) H(s,u) du \,, \quad H(t,t)=1\,,
$$
for $H(s,t)=R(s,t) \exp(\int_t^s f'(K(u))du)$,
whose unique solution is 
the Laplace transform of the semi-circle probability measure, 
evaluated at $\b(s-t)$. 
Plugging this expression in
\req{eqC} and \req{eqZ}, upon scaling both the function $f$ and
time by a factor $\b$, we recover the limiting
equation of \cite[(2.16)]{BDG} after
some integrations by parts. Further, setting $K(s)=1$ and
$\partial_s K(s)=0$ in \req{eqZ}, while replacing 
$f'(K(s))$ in \req{eqR}--\req{eqZ} by a time varying 
constant $z(s)$, corresponds to the hard
spherical constraint of \cite{CK}. Indeed, the limiting
equations of \cite{CK} are thus recovered.

Note that the constant $\beta$ can be embedded into 
$\{a_p\}$ resulting with 
$\beta G^j(\cdot) \mapsto G^j(\cdot)$
and then having $\beta=1$ in the stochastic differential 
system \req{interaction}. 
Adopting
this convention, we thus take hereafter $\beta=1$. It is 
trivial to check that the explicit dependence of 
\req{eqR}--\req{eqZ} on $\beta$ is indeed as stated,
i.e., with each appearance of $\nu'(\cdot)$, $\nu''(\cdot)$
(and $\psi(\cdot)$) multiplied by $\beta^2$.

The empirical quantities $K_N(s):=C_N(s,s)$ and 
\begin{equation}\label{eq:hfdef}
A_N(s,t):=\oneN \sum_{i=1}^N G^i(\Bx_s) x^i_t\,,\qquad\qquad\qquad
F_N(s,t):=\oneN \sum_{i=1}^N G^i(\Bx_s) B^i_t\,,
\end{equation}
play a key role in the derivation of Theorem \ref{thm-macro}. Indeed,
with
\begin{equation}\label{eq:dedef}
D_N(s,t):=-
f'(\E (K_N(t))) C_N(s,t) + A_N(t,s) \,, \qquad
E_N(s,t):=-f'(\E(K_N(s))) \chi_N(s,t) + F_N(s,t) \,,
\end{equation}
the key step of the proof of Theorem \ref{thm-macro} is summarized by  
\begin{prop}\label{prop-macro}
Fixing any $T<\infty$, in case $\b=1$, any limit point of the
sequence $(\E C_N,\E \chi_N, \E D_N, \E E_N)$ with respect to
uniform convergence on $[0,T]^2$, satisfies the integral
equations
\begin{eqnarray}
C(s,t)&=&C(s,0)+\chi(s,t)+\int_0^t D(s,u)du,\label{eqC1}\\
\chi(s,t)&=& s\wedge t+ \int_0^s E(u,t)du,\label{eq:chi}\\
D(s,t)&=&-f'(C(t,t))C(t,s)-\int_0^{t \vee s} \nu'(C(t, u))D(s,u) du
-\int_0^{t \vee s} C(s,u)\nu''(C(t,u)) D(t,u) du
\nonumber\\ &&
+C(s,t \vee s)\nu'(C(t \vee s,t))-C(s,0)\nu'(C(0,t)),
\label{eqD}\\
E(s,t)
&=&-f'(C(s,s))\chi(s,t)-\int_0^s  \nu'(C(s, u)) E(u,t) du
-\int_0^s\chi(u,t)\nu''(C(s, u))D(s,u) du\nonumber\\
&&+
\chi(s,t) \nu'(C(s,s))
-\int_0^{t \wedge s} \nu'(C(s,u)) du,
\label{eqE}
\end{eqnarray}
in the space of bounded continuous functions on $[0,T]^2$
subject to the symmetry condition $C(s,t)=C(t,s)$ and the
boundary conditions $E(s,0)=0$ for all $s$, and
$E(s,t)=E(s,s)$ for all $t \geq s$.
\end{prop}

We next detail the organization of the paper, and hence,
that of the proof of Theorem \ref{thm-macro}.

In Section \ref{sec-prelim} we prove the existence of 
strong solutions $\Bx_t$ for the Stochastic Differential System 
(in short {\bf SDS}) given in \req{interaction}, for any $N<\infty$
(see Proposition \ref{existN}). We then prove that the
functions $A_N$, $F_N$, $\chi_N$ and $C_N$
associated with these solutions have uniformly bounded (in $N$)
finite moments of all order, and
form pre-compact sequences with respect to 
uniform convergence on compacts both almost surely and in the
mean (see Proposition \ref{tightchic}). 
Applying the ``localized
concentration of measure'' of Lemma \ref{loc-conc}, we complete the
preliminary analysis of the finite size {\bf SDS} by proving
in Proposition \ref{self-average} that as $N \to \infty$ 
each of the four functions $A_N$, $F_N$, $\chi_N$ and $C_N$
``self-averages'', namely, concentrates around its mean.
These results rely on the bounding in Appendix \ref{sec-gauss}
norms of the disorder $\BJ$ that are expressed as the supremum of 
certain Gaussian fields.

The proof of Proposition \ref{prop-macro}, namely, 
that each limit point of $(\E C_N,\E \chi_N,\E D_N, \E E_N)$ must
satisfy \req{eqC1}--\req{eqE}, 
is the subject of Section \ref{sec-der}. 
Of these equations, 
upon multiplying the integrated form 
\req{eq:integ} of our {\bf SDS} by $x_t^i$ or by $B_t^i$, 
then averaging over $i$ and the probability space, 
both \req{eqC1} and \req{eq:chi} are  
immediate consequences of self-averaging. The crux of 
the proof is thus Proposition \ref{comp1}, where we show 
that as $N \to \infty$, both $\E A_N$ and $\E F_N$ are
well approximated by certain combination of our four 
functions, thereby leading to \req{eqD} and \req{eqE}.
The emergence of $\nu'(C(s,u))$ and $\nu''(C(s,u))$ in
the latter pair of equations, and hence in Theorem \ref{thm-macro},
is a consequence of the structure of the covariance kernel 
$k_{ts}(\Bx)=\E_{\BJ} [G^i(\Bx_t)G^j(\Bx_s)]$, obtained by
integrating over the disorder parameters $\BJ$ assuming 
their independence of the frozen path $\Bx_t$ 
(see Lemma \ref{kval}). However, the main difficulty of
the proof is the intricate dependence between $\BJ$ and 
$\{ \Bx_t , 0 \leq t \leq T \}$. Taking full advantage 
of the Gaussian law of $\BJ$ and the Brownian law of 
$\BB$, this difficulty 
is dealt with by combining the It\^o's calculus
identities of Appendix \ref{sec-ito} with Girsanov's theorem and 
the resulting Gaussian change of measure identities that are
derived in Appendix \ref{sec-alice}. This approach succeeds
in deriving a ``closed'' system of finitely many limiting 
equations thanks to the fact that apart from our self-averaging
global quantities, the 
kernel $k_{ts}$ 
is a quadratic form of
$\Bx$. Note that certain Hamiltonians other than \req{eq:hamil} 
also have such a property, hence are amenable to a similar treatment
(one such example is $H^N_{\BJ}(\Bx/\|\Bx\|)$ for 
$H^N_{\BJ}(\cdot)$ of \req{eq:hamil}).

Section \ref{sec-uniq} mostly deals with analytic considerations.
Its starting point is Lemma \ref{lem-diff}, showing
that any solution of \req{eqC1}--\req{eqE} is sufficiently differentiable 
to give rise to a solution of \req{eqR}--\req{eqZ}. This
is followed by Proposition \ref{uniqueness}, establishing
the uniqueness of the latter system of equations
by a Gronwall type argument. Using these two ingredients,
as well as the pre-compactness and self-averaging of our
four functions, we conclude by deducing Theorem \ref{thm-macro} out of
Proposition \ref{prop-macro}.

\section{Strong solutions, self-averaging and compactness}\label{sec-prelim}

We start with the almost sure existence of the strong solution $\Bx_t$
of (\ref{interaction}).

\begin{prop}\label{existN} 
Assume that $f'$ is locally Lipschitz, satisfying (\ref{eq:fcond}).
Then, for any $N\in\N$,  almost any $\BJ$, initial condition $\Bx_0$
and Brownian path $\BB$,
there exists a unique strong solution to (\ref{interaction}).
This solution is also unique in law for almost any $\BJ$, and $\Bx_0$,
it is a \pm on $\Ca(\R^+,\R^N)$ which we denote $\P^N_{\Bx_0, \bJ}$.
Further, with
\begin{equation}\label{jnorm}
||\bJ||_\infty^{N}=\max_{1\le p\le m}\sup_{||{\bf u^i}||\leq 
1, 1\le i\le p}\Big|
\sqrt{N}^{-1}\sum_{1\le i_k\le N, 1\le k\le p}N^{\frac{p-1}{2}}
 J_{i_1\cdots i_p}
u_{i_1}^1  \cdots u_{i_p}^p \Big| \,
\end{equation}
we have for $\d>0$ of (\ref{eq:fcond}), $q:=m/(2\d)+1$,
some $\kappa<\infty$, all $N$, $z>0$, $\bJ$, and $\Bx_0$, that
\begin{equation}\label{eq:conc}
\P^N_{{\Bx}_0,\bJ}\Big(\sup_{t\in\R^+}
K_N(t) \ge K_N(0) + \kappa (1 + \|\bJ\|_\infty^N )^q + z \Big)
\le  e^{-z N} \,.
\end{equation}
Consequently, for any $L>0$, there exists $z=z(L)<\infty$
such that
\begin{equation}\label{eq:ktail}
\P\Big(\sup_{t\in\R^+} K_N(t) \ge z \Big) \le  e^{-L N} \,.
\end{equation}
\end{prop}
\prf For every $M>0$ we introduce a 
bounded globally Lipschitz function $\phi_M$ on $\R^N$
which we choose such that $\phi_M(\bx)=\bx$
when $||\bx||\le \sqrt{NM}$, and 
then consider the truncated drift $b^M({\bf u})=(b^M_1(
{\bf u}),\ldots,b^M_N({\bf u}))$ given by 
$b^M_i({\bf u})=G^i({\bf \phi_M(u)})-f'(N^{-1}|{\bf u}|^2 \wedge M)u^i$.

Since $f'$ is  locally Lipschitz,
and since $||\bJ||_\infty^{N}$ is finite almost surely
for all $p$ and $N$ ,
it is thus clear that the drift $b^M({\bf u})$ is globally Lipschitz.
The existence and uniqueness of a square-integrable
strong solution ${\bf u}^{(M)}$ for the {\bf SDS}
$$
du^i_{t}= b^M_i ({\bf u}_t) dt  + dB^i_t
$$
is thus standard (for example, see
\cite[Theorems 5.2.5, 5.2.9]{Karatzas-shreve}).
With ${\bf u}^{(M)}$ defined for all $M$ 
on the same probability space and
filtration, consider the stopping times
$\tau_M= \inf \{ t : ||{\bf u}^{(M)}_t|| \geq \sqrt{N M} \}$. Note that
${\bf u}^{(M)}$ is the unique strong solution of (\ref{interaction}) for
$t \in [0,\tau_M]$, with $\tau_M$ a non-decreasing sequence.
By the Borel-Cantelli lemma, it suffices for the existence of a
unique strong solution ${\bf u} = \lim_{M \to \infty} {\bf u}^{(M)}$ of
the {\bf SDS} (\ref{interaction}) in $[0,T]$, to show that
\begin{equation}\label{al1}
\sum_{M=1}^\infty \P\left(\tau_M\le T\right)<\infty.
\end{equation}
To this end, fix $M$ and let 
${\bf x}_t={\bf u}^{(M)}_{t \wedge \tau_M}$
and  $Z_s=2 N^{-1} \sum_{i=1}^N \int_0^{s \wedge \tau_M} x^i_t dB^i_t$.
Applying Ito's formula for $C_N(t)=N^{-1}||{\bf x}_t||^2$ we see that
\begin{equation}
C_N(s) \leq C_N(0) + 2 
\sum_{p=1}^m \frac{a_p ||\bJ||_\infty^{N}}{(p-1)!}
\int_0^{s \wedge \tau_M} C_N(t)^{\frac{p}{2}}
 dt
-2 \int_0^{s \wedge \tau_M} f'(C_N(t)) C_N(t) dt + Z_s
+ s \wedge \tau_M  \,.
\label{al2}
\end{equation}
Since $x^{1 -{\frac{m}{2}}} f'(x) \to \infty$, it follows 
from (\ref{al2}) that there is an almost surely finite constant 
$c(||\bJ||_\infty^N)$, independent of $M$, such that 
\begin{equation}
C_N(s)
\leq C_N(0) + c(||\bJ||_\infty^N )s + Z_s
\label{al3}
\end{equation}
As the quadratic variation of the martingale $Z_s$ is
$(4/N) \int_0^{s \wedge \tau_M} C_N(t) dt \leq 4 s N^{-1} M$,
applying Doob's inequality (c.f. \cite[Theorem 3.8, p. 13]{Karatzas-shreve}) 
for the exponential 
martingale 
$L_s^\lambda 
=\exp(\lambda Z_s- 2 (\lambda^2/N) \int_0^{s \wedge \tau_M} C_N(t) dt)$
(with respect to the filtration $\{\Ha_t\}$ of $\BB_t$),
yields that 
\begin{equation}\label{Doob}
\P\left( \sup_{s\le T} \{Z_s - 2 \int_0^s C_N(t) dt\}\ge z \right)
\le\P\left( \sup_{s\le T} L_s^N \geq e^{z N} \right) \le e^{-z N}
\,,
\end{equation}
for any $z>0$. Therefore, (\ref{al3})
shows that with probability greater than $1- e^{-z N}$, 
$$
C_N(s \wedge \tau_M) \le C_N(0) + c(||\bJ||_\infty^N ) T + z
 +2 \int_0^{s \wedge \tau_M} C_N(t) dt \;,
$$
for all $s\le T$,
and by Gronwall's lemma then also
\begin{equation}
\sup_{t\le T}N^{-1} |{\bf u}^{(M)}_{t \wedge \tau_M} |^2
\le [C_N(0) + c(||\bJ||_\infty^N ) T + z]e^{2 T} .
\label{new}
\end{equation}
Setting $z=M/3$,
for large enough $M$ (depending of $N$, $\bJ$, $\Bx_0$
and $T$ which are fixed here), 
the right-side of (\ref{new}) is at most 
$M/2$, resulting with 
$$
\P(\tau_M\le T) \leq e^{-M N/3},
$$
and hence with (\ref{al1}).
We also have  weak uniqueness
of our solutions for almost
all ${\bf J}$ since the restriction of
any weak solution to the stopped $\sigma$-field $\Ha_{\tau_M}$
for the filtration $\Ha_t$ of $\BB_t$ is unique. We denote 
this unique  weak solution of (\ref{interaction}) 
by $\P^{N}_{\Bx_0,\BJ}$.

Turning to the proof of (\ref{eq:conc}),
by (\ref{eq:fcond}), for any $c>0$
there exists $\kappa<\infty$ such that for all $r,x \geq 0$,
$$2[ f'(x)x - r \sum_{p=1}^m \frac{a_p x^{\frac{p}{2}}}{(p-1)!}] -1
\ge c x- \kappa (1+r)^q.$$
Taking $r= \|\bJ\|_\infty^{N}$,
we see that by (\ref{al2}), for all $N$ and $s \geq 0$,
$$
C_N(s\wedge\tau_M) \leq C_N(0) -
\int_0^{s \wedge \tau_M}  [ c C_N(t) - \kappa (1+ \|\bJ\|_\infty^{N})^q ] dt
 + Z_s \,,
$$
where $(Z_s)_{s\ge 0}$ is a martingale with bracket
$( 4 N^{-1}\int_0^{s\wedge \tau_M} C_N(t) dt,s\ge 0)$.

By Doob's inequality (\ref{Doob}), with probability
at least $1-e^{- z N}$, 
$$\sup_{u\le s \wedge \tau_M} Z_u \le 2 \int_0^{s \wedge \tau_M}
C_N(t)dt+z,$$
for all $s \geq 0$. Setting $c=3$ we then have that
\begin{eqnarray}
C_N(s\wedge\tau_M) &\leq& C_N(0) +z 
-\int_0^{s \wedge \tau_M}  C_N(t) dt 
+ \kappa  (1+ \|\bJ\|_\infty^{N})^q (s \wedge \tau_M) \,,
\label{al2bb}
\end{eqnarray}
so that by Gronwall's lemma,
$$C_N(s\wedge\tau_M)\le e^{-s \wedge \tau_M} (C_N(0)+z)
+ \kappa  (1+ \|\bJ\|_\infty^{N})^q 
\int_0^{s \wedge \tau_M} e^{-t} dt$$
from which the conclusion (\ref{eq:conc}) is obtained by
considering $M \to \infty$. 

In view of the assumed
exponential in $N$ decay of the tail probabilities for $K_N(0)$
and the bound \req{bovier2} on the corresponding probabilities for
$\|\bJ\|_\infty^N$ we thus get also the bounds of \req{eq:ktail}.
\qed

Our next lemma provides bounds on $G^i(\bx)$ of (\ref{eq:gdef})
which we often use en-route to the uniform boundedness of moments, 
pre-compactness, and concentration around the mean 
of the functions $A_N$, $F_N$, $\chi_N$, and $C_N$.
\begin{lem}\label{Ginc}
There exists a constant $c<\infty$ such that, for all
$N<\infty$ and every $\bx$, $\wbx \in \reals^N$, 
\begin{equation}
\label{eq:gbd1}
\Big[ \sum_{i=1}^N (G^i(\bx)-G^i(\wbx))^2 \Big]^{\frac{1}{2}}\le
c||\bJ||_\infty^N [1+({\frac{1}{N}}||\bx||^2)^{\frac{m-2}{2}}] 
[1+({\frac{1}{N}}||\wbx||^2)^{\frac{m-2}{2}}] ||\bx-\wbx|| \,,
\end{equation}
and in particular, 
\begin{equation}
\label{eq:gbd2}
\Big[{\frac{1}{N}}\sum_{i=1}^N G^i(\bx)^2 \Big]^{\frac{1}{2}}\le
c||\bJ||_\infty^N [1+({\frac{1}{N}}||\bx||^2)^{\frac{m-1}{2}}] \,.
\end{equation}
\end{lem}
\prf Fixing $N$, $\bx$, $\wbx$ and $\by \in \reals^N$,
note the telescoping sum
$$
\sum_{i=1}^N (G^i(\bx)-G^i(\wbx)) y^i
=
\sum_{p=1}^{m} \frac{a_{p}}{(p-1)!} \sum_{l=1}^{p-1}
\sum_{i_k}J_{i_1\cdots i_p}(\prod_{h=1}^{l-1} x^{i_h})(x^{i_{l}}
-\tilde x^{i_l})( \prod_{h=l+1}^{p-1} \tilde  x^{i_{h}}) y^{i_p}
$$
Recall the definition (\ref{jnorm}) leading to the bounds,
$$
| \sum_{i_k}J_{i_1\cdots i_p}(\prod_{h=1}^{l-1} x^{i_h})(x_{i_{l}}
-\tilde x_{i_l})(\prod_{h=l+1}^{p-1} \tilde  x^{i_{h}}) y^{i_p}|\le
||\bJ||_\infty^N ({\frac{1}{ N}}||\bx||^2)^{\frac{l-1}{2}} 
({\frac{1}{ N}}||\wbx||^2)^{\frac{p-l-1}{2}} ||\bx-\wbx|| ||\by|| \,,
$$
for $1 \leq l \leq p-1 \leq m-1$.
Consequently, we get for some $c=c(a_1,\ldots,a_m)<\infty$ which
is independent of $N$, that 
$$
\sum_{i=1}^N (G^i(\bx)-G^i(\wbx)) y^i \le c ||\bJ||_\infty^N 
[1+({\frac{1}{N}}||\bx||^2)^{\frac{m-2}{2}}] 
[1+({\frac{1}{N}}||\wbx||^2)^{\frac{m-2}{2}}]
||\bx-\wbx|| ||\by|| \,,
$$
for all $\bx$, $\wbx$ and $\by$. Taking $y^i=G^i(\bx)-G^i(\wbx)$ results
with the bound of (\ref{eq:gbd1}), out of which we get 
the bound (\ref{eq:gbd2}) by considering $\wbx={\bf 0}$.
\qed

By the estimate (\ref{bovier1}) of Appendix 
\ref{sec-gauss}, we know that for any $k<\infty$,
\begin{equation}\label{eq:ger2}
\sup_N \E [ (||\bJ||_\infty^{N})^k ] < \infty \,,
\end{equation}
for the norm $||\bJ||_\infty^{N}$ of (\ref{jnorm}).
By the boundedness of $N \mapsto \E [ K_N(0)^k ]$ and
(\ref{eq:conc}) the bound \req{eq:ger2} immediately implies that 
for each $k<\infty$, also
\begin{equation}\label{eq:conc2}
\sup_{N} \E \Big[ \sup_{t\in\R^+} K_N(t)^{k} \Big] <\infty\,.
\end{equation}

Building upon \req{eq:ger2} and \req{eq:conc2} we
next prove uniform moment bounds and pre-compactness of the 
 functions of interest to us.
\begin{prop}\label{tightchic}
Let $U_N$ denote any one of the functions 
$A_N$, $F_N$, $\chi_N$ and $C_N$. Then, 
$\E \big[ \sup_{s,t \leq T} |U_N(s,t)|^k \big]$ is uniformly bounded in
$N$, for any fixed $T<\infty$ and $k<\infty$. Further, 
for any fixed $T<\infty$, the sequence of 
continuous functions $ U_N(s,t)$ is pre-compact 
almost surely and in expectation
with respect to the uniform topology on $[0,T]^2$. 
\end{prop}
\prf We start with the uniform bound on the moments of 
$A_N$, $F_N$, $\chi_N$ and $C_N$. To this end, let 
$B_N(t):=\oneN \sum_{i=1}^N |B_t^i|^2$ and 
$G_N(t):=\oneN \sum_{i=1}^N |G^i(\Bx_t)|^2$.
Fixing
$T,k<\infty$ let $\|K_N\|_\infty:=\sup \{ K_N(t) : 0 \leq t \leq T \}$, and
similarly define $\|B_N\|_\infty$, $\|G_N\|_\infty$ and
$\|U_N\|_\infty:=\sup \{ U_N(s,t) : 0 \leq s,t \leq T \}$. 

A key estimate is then the following,
\begin{equation}\label{ger:moments} 
\sup_N \E [ (||\bJ||_\infty^{N})^k ] +
\sup_N\E[\|K_N\|_\infty^k] 
+\sup_N\E[\|B_N\|_\infty^k] + \sup_N\E[\|G_N\|_\infty^k]<\infty \,,
\end{equation}
holding for each fixed $k$. Indeed,
the bounds on $\|\bJ\|_\infty^N$ and
$\|K_N\|_\infty$ are already obtained in \req{eq:ger2} and
\req{eq:conc2},
whereas by \req{eq:gbd2} 
we have that
\begin{equation}\label{eq:gubd} 
(G_N(t))^{\frac{1}{2}}\leq c ||\bJ||_\infty^N [1+ K_N(t)^{\frac{m-1}{2}}]\,,
\end{equation}
yielding by \req{eq:ger2} and \req{eq:conc2}
the 
uniform moment bound on $\|G_N\|_\infty$. Finally, 
it is easy to show that 
\begin{equation}\label{eq:bbd}
\P (\|B_N\|_\infty \geq 4T(3+x))\leq e^{-N x}
\,,
\end{equation}
for all 
$T,N<\infty$ and $x>0$ (c.f. \cite[(3.44)]{BDG}), 
thereby providing a uniform bound for each moment of $\|B_N\|_\infty$
and concluding the derivation of \req{ger:moments}.

Similarly, by \req{eq:ktail}, \req{eq:gubd},
\req{eq:bbd} and the 
exponential in $N$ bound of \req{bovier2} on the tail of
$\|\bJ\|^N_\infty$, we have for each $L> 0$, that 
there exists $M=M(L) < \infty$ such that for all $N$,
\begin{equation}\label{ger:moments2} 
\P (\|\bJ\|^N_\infty + \|K_N\|_\infty 
+\|B_N\|_\infty + \|G_N\|_\infty\ge M) \le e^{-LN} \,.
\end{equation}

In view of \req{ger:moments} and \req{ger:moments2},
we get the uniform in $N$ bounds on moments of 
$\|U_N\|_\infty$ and the exponential in $N$ bounds on the tail 
probabilities for $\|U_N\|_\infty$ upon observing that 
\begin{eqnarray*}
|C_N(s,t)|   &\leq& K_N(s)+K_N(t), \quad 
|\chi_N(s,t)| \leq K_N(s)+B_N(t), \quad \\
|A_N(s,t)|   &\leq& G_N(s)+K_N(t), \quad 
|F_N(s,t)| \leq G_N(s)+B_N(t). 
\end{eqnarray*}

With the previous controls
on  $ \| U_N\|_\infty$ already 
established, by the Arzela-Ascoli theorem, the 
pre-compactness of $ U_N$ follows by showing that it is 
an equi-continuous sequence. To this end, observe that such
$U_N(s,t)$ are all of the form $\oneN \sum_{i=1}^N a^i_s b^i_t$
hence, 
\begin{eqnarray*}
&& 
|  U_N(s,t) -  U_N(s',t') |\leq
  \oneN \sum_{i=1}^N  |a^i_s-a^i_{s'}| |b^i_t| 
+  \oneN \sum_{i=1}^N |a^i_{s'}| |b^i_t-b^i_{t'}| \\
&\leq& 
\big[\oneN \sum_{i=1}^N  |a^i_s-a^i_{s'}|^2 \big]^{1/2}
\big[\oneN \sum_{i=1}^N  |b^i_t|^2\big]^{1/2}
+ 
\big[\oneN \sum_{i=1}^N  |b^i_t-b^i_{t'}|^2\big]^{1/2} 
\big[\oneN \sum_{i=1}^N  |a^i_{s'}|^2\big]^{1/2} \,.
\end{eqnarray*}
Here the functions $\ba_t$ and $\bb_t$ are either $\Bx_t$,
$\BB_t$ or $G(\Bx_t)$. So, in view of \req{ger:moments}
and \req{ger:moments2}, it suffices to show that for 
any $\e>0$, some function $L(\d,\e)$ going 
to infinity as $\d$ goes to zero and all $N$, 
$$
\P(\sup_{|t-t'|<\delta}  \Big[ \oneN \sum_{i=1}^N |b^i_t-b^i_{t'}|^2 
\Big] >\e)\le e^{-L(\d,\e)N} , \quad
\sup_{|t-t'|<\delta}  \E\Big[ \oneN \sum_{i=1}^N |b^i_t-b^i_{t'}|^2 
\Big]\le L(\d,\e)^{-1}\,,
$$
where $\bb=\Bx$, $\BB$ or $G(\Bx)$. Obviously, this holds for 
$\bb=\BB$, in which case the expectation equals $\delta$
regardless of $N$ and the tail probability bound follows 
upon considering the union of such probabilities 
for $t,t' \in [i\delta,(i+2)\delta]$,
$i=0,1,\ldots,T/\delta$ and applying 
\req{eq:bbd} with $T=2\delta$.
In case $\bb=G(\Bx)$, note that by \req{eq:gbd1} we have that
$$
\frac{1}{N}\sum_{i=1}^N (G^i(\bx_t)-G^i(\bx_{t'}))^2 \leq 4c^2
(\|\BJ\|^{N}_\infty)^2 (\frac{1}{N}\|\bx_t-\bx_{t'}\|^2) 
(1+K_N(t)^{(m-2)})(1+K_N(t')^{(m-2)}) 
$$
for all $t,t'$. Thus, in view of the bounds
\req{ger:moments} and \req{ger:moments2},
everything reduces to showing that
\begin{equation}\label{eq:only}
\P( \sup_{|t-t'|<\delta} \oneN \sum_{i=1}^N |x^i_t-x^i_{t'}|^2>\e)
\le e^{-L'(\d,\e) N}, 
\sup_{|t-t'|<\delta} \E \Big[ \big( 
\oneN \sum_{i=1}^N |x^i_t-x^i_{t'}|^2 \big)^2 \Big] \le L'(\d,\e)^{-1}
\end{equation}
for some $L'$ with the same properties as $L$. To this end,
 note that by \req{interaction}
$$|x^i_t-x^i_{t'}|\le 
 |B^i_t-B^i_{t'}| + \|f'(K_N)\|_\infty
\int_t^{t'}|x^i_u|du 
+\int_t^{t'}|G^i (\Bx_u)|du \,.
$$
So by \req{eq:fcondu}, 
for some universal constant $\rho_1<\infty$, all $t,t'$ and $N$,
$$
\oneN \sum_{i=1}^N |x^i_t-x^i_{t'}|^2 \leq 
\frac{3}{N} \sum_{i=1}^N |B^i_t-B^i_{t'}|^2 
+3 |t-t'|^2 \Big[ \rho_1 (1+\|K_N\|_\infty)^{2r}\|K_N\|_\infty
+\|G_N\|_\infty \Big]
$$
Thus, by \req{ger:moments}, \req{eq:bbd} and \req{ger:moments2},
we readily obtain \req{eq:only}.
\hfill\qed

A key ingredient in the derivation of the limiting dynamics is 
the ``self-averaging'' of the functions $A_N$, $F_N$, $C_N$, 
$\chi_N$ (and hence of $D_N$ and $E_N$), as we next state and prove. 
\begin{prop}\label{self-average}
Assume that the concentration of measure of Hypothesis 
\ref{defconcentration} holds,
and as before let $U_N$ denote any one 
of the functions $A_N$, $F_N$, $\chi_N$ and $C_N$. Then, 
for any $T<\infty$ and $\rho > 0$,
\begin{equation}\label{eq:asself}
\sum_N \P [ \sup_{s,t \le T} |U_N(s,t)- \E(U_N(s,t))| \geq \rho] < \infty \,,
\end{equation}
implying by the uniform moment bounds on $\|U_N\|_\infty$ that 
\begin{equation}\label{eq:l2self}
 \lim_{N \rightarrow \infty} 
\sup_{s,t \leq T} \E \Big[ |U_N(s,t)-\E U_N(s,t)|^{2} \Big]=0\,.
\end{equation}
\end{prop}

Our strategy relies on the following general 
"localized" version of the concentration of measure property, which 
might be of some independent interest.
\begin{lem}\label{loc-conc}
Suppose $V_N$ are functions on normed spaces 
$\Ea_N$ in which the concentration of measure
of Hypothesis \ref{defconcentration} holds. 
 Assume that 
$K=\sup_N \E[V_N^2]^{1/2}<\infty$ and 
that for every $M<\infty$ there exists a subset 
$\CL_{N,M}$ of $\Ea_N$ on which $|V_N| \leq 2M$
and $V_N$ is Lipschitz with Lipschitz constant 
$A_{N,M} \leq \frac{D(M)}{\sqrt{N}}$.  Further 
assume that, for every $L>0$ there exists an $M=M(L)$ 
such that 
$\P[\CL_{N,M}^c]\le \exp(-LN)$ for all $N$. Then,
\begin{equation}\label{concentration}
\P[ |V_N-\E[V_N]|\geq \rho ] \leq C^{-1}
\exp{(- C({\frac{\rho}{2 D(M(L))}})^\alpha)
N^{\frac{\alpha}{2}})} + 4 (K+M(L)) \rho^{-1} e^{-LN/2} + e^{-NL} \;.
\end{equation}
\end{lem}
\prf 
Recall that for every $M<\infty$ there exists a globally 
Lipschitz function 
$$
U_N^M(z)=\sup_{y \in \CL_{N,M}} \{V_N(y)-A_{N,M} \|z-y\|\}
$$ 
on $\Ea_N$, of Lipschitz constant $A_{N,M}$,
which coincides with $V_N$ on the set $\CL_{N,M}$.
These properties are inherited by $V_N^M = \max(U_N^M,-2M)$
for which also $|V_N^M| \le 2M$.
We thus get \req{concentration} 
by combining
the triangle inequality
$$ 
|V_N-\E[V_N]| \le  |V_N-V_N^M| + |V_N^M-\E[V_N^M]| 
+ |\E[V_N^M]-\E[V_N]| \,, 
$$
with the Cauchy-Schwartz inequality
$$ |\E[V_N-V_N^M]| \leq \E |V_N-V_N^M| \leq  2M \P(\CL_{N,M}^c) +
 \E[|V_N| {\bf 1}_{\CL_{N,M}^c} ]
\leq (K+2M) \P[\CL_{N,M}^{c}]^{1/2}\,,$$
and applying the concentration of measure inequality of Hypothesis 
\ref{defconcentration} for $V_N^M$. 
\qed

\noindent
{\bf Proof of Proposition \ref{self-average}:}
We wish to apply the estimate \req{concentration} to $V_N=U_N(s,t)$ 
for any of our four functions, and any fixed pair of times $s,t$.
To this end, for each $M<\infty$ and any $N$ define the subset  
$$
\CL_{N,M} = \{
(\Bx_0, \BJ, \BB)\in \Ea_N \,:\, 
||\BJ||^{N}_\infty + \|B_N\|_\infty + \| K_N\|_\infty + \| G_N \|_\infty
\le M \;\} 
$$
of $\Ea_N$. 
For $M$ sufficiently large, the probability of the complement set
$\CL_{N,M}^c$ decays exponentially in $N$ by \req{ger:moments2},
whereas by Proposition \ref{tightchic}
we have the uniform moment bounds for the functions 
$U_N(s,t)$, as well as the stated pointwise bound in $\CL_{N,M}$. 
It thus suffices to prove the stated Lipschitz constant 
of $V_N$ on $\CL_{N,M}$.  We start by showing next that restricted to 
the set $\CL_{N,M}$, the solution $\bx$ of (\ref{interaction}) is
a Lipschitz function of the triple $(\Bx_0, \BJ, \BB)$.
\begin{lem}\label{x is Lipschitz} 
Let $\bx,\widetilde{\bx}$ be the two strong solutions
of (\ref{interaction}) constructed from $(\Bx_0,\BJ,\BB)$
and $(\widetilde{\Bx}_0,\widetilde{\BJ}, \widetilde{\BB})$, respectively.
If $(\Bx_0,\BJ,\BB)$
and $(\widetilde{\Bx}_0,\widetilde{\BJ},\widetilde{\BB})$ 
are both in $\CL_{N,M}$, then
we have the Lipschitz estimate:
\begin{equation}\label{eq:lipbd}
\sup_{t\le T}\frac{1}{N}\sum_{1\le i\le N}
|x^i_t-\widetilde{x}^i_t|^2\le \frac{D_o(M,T)}{N}\|(\Bx_0, \BJ, \BB) 
- (\widetilde{\Bx}_0,\widetilde{\BJ},\widetilde{\BB})\|^2\,,
\end{equation}
where the finite 
constant $D_o(M,T)$ is independent of $N$.
\end{lem}
\prf
We denote by $G^i(\cdot)$ (resp. $\tilde G^i(\cdot)$)
the Gaussian fields constructed from the $\BJ$ (resp. $\tilde\BJ$).
We write the following natural decomposition:
\begin{eqnarray}
\label{pluc}
e_N(t):= {\frac{1}{ N}}\sum_{i=1}^N
|x^i_t-\widetilde{x}^i_t|^2
&=&{\frac{1}{ N}}\sum_{i=1}^N
(x^i_t-\widetilde{x}^i_t)\bigg( 
 (x^i_0-\widetilde{x}^i_0)+(B^i_t-\widetilde{B}^i_t)
\nonumber
\\
&&- \int_0^t (f'(K_N(u))-f'(\widetilde K_N(u))) x^i_udu
-\int_0^t f'(\widetilde K_N(u)) ( x^i_u-\widetilde x^i_u)du
\\
&&+ \int_0^t (G^i(\bx_u)-G^i(\widetilde \bx_u)) du
+\int_0^t (G^i(\widetilde \bx_u)-\widetilde G^i(\widetilde \bx_u)) du
\bigg)\nonumber\\
&=& I_1+I_2+\cdots+I_6
\nonumber
\end{eqnarray}
Then, since $\| K_N\|_\infty \le M$ and $\|\tilde K_N\|_\infty\le M$,
we have that for some finite $C(M,T)$ that is 
independent of $N$ and for all $c>0$,
\begin{eqnarray*}
I_1&\le& {\frac{c}{2N}}\sum_{i=1}^N
(x^i_t-\widetilde{x}^i_t)^2 +{\frac{1}{2cN}}\sum_{i=1}^N
(x^i_0-\widetilde{x}^i_0)^2\\
I_2&\le&{\frac{c}{2N}}\sum_{i=1}^N
(x^i_t-\widetilde{x}^i_t)^2 +{\frac{1}{2cN}}\sum_{i=1}^N
(B^i_t-\widetilde{B}^i_t)^2\\
I_3+I_4&\le& C(M,T)\Big({\frac{c}{2N}}\sum_{i=1}^N
(x^i_t-\widetilde{x}^i_t)^2 +\int_0^t {\frac{1}{2cN}}\sum_{i=1}^N
(x^i_u-\widetilde{x}^i_u)^2 du \Big)\\
I_6&\le& 
C(M,T)\Big({\frac{c}{2N}}\sum_{i=1}^N
(x^i_t-\widetilde{x}^i_t)^2+{\frac{1}{2cN}}\sum_{p=1}^m \sum_{1 \leq i_1\ldots i_{p} \leq N} 
(N^{\frac{p-1}{2}}(J_{i_1\cdots i_{p}}-\tilde J_{i_1\cdots i_{p}}))^2 \Big)
\end{eqnarray*}

We can bound the term $I_5$ on $\CL_{N,M}$, using \req{eq:gbd1} of
Lemma \ref{Ginc} by
$$I_5\le C(M)||\bJ||_\infty^N \int_0^t  \left({\frac{1}{ N}}\sum_{i=1}^N
(x^i_t-\widetilde{x}^i_t)^2 \right)^{\frac{1}{2}}\left({\frac{1}{ N}}\sum_{i=1}^N
(x^i_u-\widetilde{x}^i_u)^2 \right)^{\frac{1}{2}}du$$
$$\le C(M,T)\Big({\frac{c}{2N}}\sum_{i=1}^N
(x^i_t-\widetilde{x}^i_t)^2 + \int_0^t 
{\frac{1}{2cN}}\sum_{i=1}^N
(x^i_u-\widetilde{x}^i_u)^2 du \Big)\,,$$
for some $C(M,T)$ independent of $N$ and all $c>0$. 
Adding these estimates we get the bound
$$
e_N(t) \leq \widetilde{C}(M,T) \Big[ c e_N(t) + \frac{1}{cN} 
\| (\bx_0,\bJ,\BB)  
- (\widetilde{\Bx}_0,\widetilde{\BJ},\widetilde{\BB})\|^2
+ \frac{1}{c} \int_0^t e_N(u) du \Big]\,,
$$
on $e_N(t)$ of \req{pluc}, so for $c=c(M,T)>0$ small enough,
Gronwall's lemma yields the stated bound \req{eq:lipbd}.
\qed

Equipped with Lemma \ref{x is Lipschitz} it is now easy to prove that 
\begin{lem}\label{U are Lipschitz} 
Let $\bx,\widetilde{\bx}$ be the two strong solutions
of (\ref{interaction}) constructed from $(\Bx_0,\BJ,\BB)$
and $(\widetilde{\Bx}_0,\widetilde{\BJ}, \widetilde{\BB})$, respectively.
If $(\Bx_0,\BJ,\BB)$
and $(\widetilde{\Bx}_0,\widetilde{\BJ},\widetilde{\BB})$ are both in 
$\CL_{N,M}$, then
we have the Lipschitz estimate for each of the four functions $U_N(s,t)$
of interest, 
\begin{equation}\label{eq:lippr}
\sup_{s,t \le T}|U_N(s,t)- \widetilde{U_N}(s,t)| 
\le \frac{D(M,T)}{\sqrt{N}} 
\|(\Bx_0, \BJ, \BB) - (\widetilde{\Bx}_0,\widetilde{\BJ},\widetilde{\BB})\|\,,
\end{equation}
where the constant $D(M,T)$ depends only on $M$ and $T$ and not on $N$.
\end{lem}

\prf
Since each of the four functions
$U_N(s,t)$ is of the form $\oneN \sum_{i=1}^N a^i_s b^i_t$, we have that
\begin{eqnarray*}
&& 
| U_N(s,t) - \widetilde{U_N}(s,t) |\leq
\oneN \sum_{i=1}^N  |a^i_s-\tilde{a}^i_{s}| |b^i_t| 
+ \oneN \sum_{i=1}^N |\tilde{a}^i_{s}| |b^i_t-\tilde{b}^i_{t}| \\
&\leq& 
\big[\oneN \sum_{i=1}^N  |a^i_s-\tilde{a}^i_{s}|^2 \big]^{1/2}
\big[\oneN \sum_{i=1}^N  |b^i_t|^2\big]^{1/2}
+ 
\big[\oneN \sum_{i=1}^N  |b^i_t-\tilde{b}^i_{t}|^2\big]^{1/2} 
\big[\oneN \sum_{i=1}^N  |a^i_{s}|^2\big]^{1/2} \,.
\end{eqnarray*}
Here the functions $\ba_t$ and $\bb_t$ are either $\Bx_t$,
$\BB_t$ or $G(\Bx_t)$. 
This bound and Lemma \ref{x is Lipschitz} are sufficient to prove the conclusion of the lemma for the two functions $C_N(s,t)$ and $\chi_N(s,t)$.
To prove it for the other two functions $A_N(s,t)$ and $F_N(s,t)$,
note that by \req{eq:gbd1} we have the bound 
$$ \big[\oneN \sum_{i=1}^N  |G^i(\Bx_s)-G^i(\tilde{\Bx}_s)|^2 \big]^{1/2} 
\le C(M)
 ||\bJ||_\infty^N (1+ M^{\frac{p-1}{2}}) \Big(\frac{1}{N}\sum_{1\le i\le N}
|x^i_s-\widetilde{x}^i_s|^2 \Big)^{1/2}\,,
$$
holding on $\CL_{N,M}$, so combining
the obvious consequence of Cauchy-Schwartz
$$\big[\oneN \sum_{i=1}^N  |G^i(\wbx_s)-\tilde{G}^i(\wbx_s)|^2 \big]^{1/2} 
\le C(M) \big[ \frac{1}{N} \sum_{p=1}^m \sum_{1 \leq i_1\ldots i_{p} \leq N} 
(N^{\frac{p-1}{2}}(J_{i_1\cdots i_{p}}
-\tilde J_{i_1\cdots i_{p}}))^2  \big]^{1/2}\,,$$
on $\CL_{N,M}$, with Lemma \ref{x is Lipschitz}, one gets
$$\big[\oneN \sum_{i=1}^N  |G^i(\Bx_s)-\tilde{G}^i(\tilde{\Bx}_s)|^2 \big]^{1/2} \le \frac{C(M,T)}{\sqrt{N}} \|(\Bx_0, \BJ, \BB) - (\widetilde{\Bx}_0,\widetilde{\BJ},\widetilde{\BB})\|\,.$$
Moreover, the estimate \req{eq:gbd2} gives the bound
$$\big[\oneN \sum_{i=1}^N  |G^i(\Bx_t)|^2\big]^{1/2} \le  
c ||\bJ||_\infty^N (1+M^{m-1}) \le C(M) \,.$$ 
The last two estimates  and Lemma \ref{x is Lipschitz} 
thus 
yield the conclusion \req{eq:lippr}
for both $A_N(s,t)$ and $F_N(s,t)$.
\qed

In view of Lemma \ref{U are Lipschitz}, the inequality  
\req{concentration} applies for $V_N=U_N(s,t)$, for any
fixed $s,t \leq T$ with constants $K$ and $D=D(M(L),T)$
that are independent of $s,t$, $\rho$ and $N$. Consequently,
by the union bound, for any finite subset $\Aa$ of $[0,T]^2$,
and any $\rho>0$, the sequence
$N \mapsto \P[ \sup_{(s,t) \in \Aa} |U_N(s,t)- \E U_N(s,t)| \geq \rho/3]$ 
is summable. 

Recall that in the course of proving Proposition \ref{tightchic}
we showed that for any $\e>0$ there exists $\widetilde{L}(\d,\e) \to \infty$ 
for $\d \to 0$, such that for all $N$,
$$
\P(\sup_{|s-s'|+|t-t'|<\d} |U_N(s,t)-U_N(s',t')|>\e) \leq e^{-
\widetilde{L}(\d,\e) N}\,, \quad
\sup_{|s-s'|+|t-t'|<\d} |\E U_N(s,t)- \E U_N(s',t')| \leq 
\widetilde{L}(\d,\e)^{-1} \,. 
$$
Choosing $\d>0$ small enough so that $\widetilde{L}(2\d,\rho/3) > 3/\rho>0$,
we thus get \req{eq:asself} by considering the finite subset 
$\Aa=\{(i\d,j\d): i,j=0,1,\ldots,T/\d\}$
of all points of $[0,T]^2$ on a $\d$-mesh.
%
\qed

We shall often apply the following direct consequence of 
Propositions \ref{tightchic} and  \ref{self-average}.
\begin{cor}\label{cor-self}
Suppose that
$\Psi: \R^\ell \to \R$ is locally Lipschitz with
$|\Psi(z)| \leq M \|z \|_k^k$ for some $M,\ell,k<\infty$,
and $Z_N \in \R^\ell$
is a random vector, where for $j=1,\ldots,\ell$,
the $j$-th coordinate of $Z_N$ is one of 
the functions $A_N$, $F_N$, $\chi_N$ or $C_N$ 
evaluated at some $(s_j,t_j) \in [0,T]^2$. Then,
$$
\lim_{N \to \infty} 
\sup_{s_j,t_j} |\E \Psi(Z_N) - \Psi(\E Z_N)| = 0 \;.
$$
\end{cor}
\prf It follows from Proposition \ref{tightchic} that
$R = \sup_{s_j,t_j,N} \|\E(Z_N) \|_k < \infty$. For each $r \geq R$
let $c_r$ denote the finite Lipschitz constant of
$\Psi(\cdot)$ (with respect to $\|\cdot\|_2$), 
on the compact set $\Gamma_r := \{ z : \|z\|_k \leq r \}$.
Then, 
\begin{eqnarray*}
|\E \Psi(Z_N) - \Psi(\E Z_N)|&\leq&
\E |\Psi(Z_N) - \Psi(\E Z_N)| {\bf 1}_{Z_N \in \Gamma_r}  + 
\E |\Psi(Z_N)| {\bf 1}_{Z_N \notin \Gamma_r} 
+|\Psi(\E Z_N)| \P (Z_N \notin \Gamma_r)  \\
&\leq& c_r \E[\|Z_N - \E Z_N \|_2] + 2 \ell M r^{-k} \E \|Z_N \|_{k}^{2k} \,.
\end{eqnarray*}
We have by \req{eq:l2self} and
the uniform moment bounds of Proposition \ref{tightchic}
that
$\sup_{s_j,t_j} \E[\|Z_N - \E Z_N \|_2] \to 0$ 
as $N \to \infty$, while 
$c' = \sup_{s_j,t_j,N} \E \|Z_N \|_{k}^{2k} < \infty$,
implying that 
$$
\lim_{N \to \infty}
\sup_{s_j,t_j} |\E \Psi(Z_N) - \Psi(\E Z_N)| \leq 2 c' \ell M r^{-k} 
\,,
$$
which we make arbitrarily small by taking $r \to \infty$.
\qed

\section{Limiting equations: proof of 
Proposition \ref{prop-macro}}\label{sec-der}

We shall denote in short
$$C_N^a(s,t)=\E[C_N(s,t)],\qquad \chi_N^a(s,t)=\E[\chi_N(s,t)]$$
where expectation is over the Brownian path $\BB$, the disorder 
$\bJ$ and the initial condition $\bx_0$, and adopt a similar
notation for other functions of interest.

Integrating the {\bf SDS} \req{interaction} we have that 
\begin{equation}\label{eq:integ}
x^i_s=x^i_0+B^i_s -\int_0^s f'(K_N(u))x^i_u du +\int_0^s
G^i(\bx_u) du \,.
\end{equation}
Hence, upon multiplying by $x^i_t$ and 
$B^i_t$ followed by averaging over $i$ and taking 
the expected value, we get that for any $s,t\in\R^+$,
\begin{eqnarray}
C_N^a(s,t)&=& C_N^a(0,t)+\chi_N^a(t,s)
-\int_0^s \E[ f'(K_N(u)) C_N(u,t)]du
+\int_0^s A_N^a (u,t) du\label{eqCa}\\
\chi_N^a(s,t)&=& \chi_N^a(0,t) +
t\wedge s-\int_0^s \E[ f'(K_N(u)) \chi_N(u,t)]du +\int_0^s F_N^a (u,t) du \,.
\label{eqchi}
\end{eqnarray}
In the following, we use $a_N\simeq b_N$ when
$a_N(\cdot,\cdot)-b_N(\cdot,\cdot) \to 0$ as $N \to \infty$,
uniformly on $[0,T]^2$.
Applying Corollary \ref{cor-self} 
(for $\Psi(z)=z_1 f'(z_2)$ whose polynomial growth is guaranteed by 
our assumption \req{eq:fcondu}), we deduce that
$$
\E[ f'(K_N(u)) C_N(u,t)]\simeq f'(K_N^a(u)) C_N^a(u,t),\quad
\E[ f'(K_N(u)) \chi_N(u,t)] \simeq f'(K_N^a(u)) \chi_N^a(u,t).
$$

Our next proposition, approximates the terms
$A_N^a$ and $F_N^a$ which we need in order 
to compute the limits
of \req{eqCa} and \req{eqchi} as $N \to \infty$.
\begin{prop}\label{comp1}
We have that 
\begin{eqnarray}\label{eq:apxa}
A_N^a(t,s)&\simeq&
\nu'(C_N^a(t,t \vee s)) C_N^a(s,t \vee s) -
\nu'(C_N^a(t,0))C_N^a(s,0) 
- \int_0^{s\vee t} \nu'(C_N^a(t,u)) D_N^a(s,u) du \nonumber \\
&-&\int_0^{s\vee t} \nu''(C_N^a(t,u))C_N^a(s,u) D_N^a (t,u) du \,.
\end{eqnarray}
Further, 
\begin{eqnarray}\label{eq:apxf}
F_N^a (s,t)&\simeq& \chi_N^a(s,t \wedge s) \nu'(C_N^a(s,s)) - 
\int_0^{t \wedge s} \nu'(C_N^a(s,u)) du 
- \int_0^s \nu'(C_N^a(s,u)) E_N^a (u,t \wedge u) du \nonumber \\
&-& \int_0^s \chi_N^a (u,t \wedge u) \nu''(C_N^a(s,u)) D_N^a(s,u) du \;.
\end{eqnarray}
\end{prop}

Deferring the proof of Proposition \ref{comp1}, we first
use its conclusion to 
complete the proof of Proposition \ref{prop-macro}. To this end, note that
$\BB_0=0$ so $F_N^a(s,0)=0$, for all $N$ and $s$. Further,
with $\Gga_s=\sigma(\BJ,\Bx_0,\BB_u, u \leq s)$, we have by
\req{eq:hfdef} that for all $N$ and  $t > s$,
$$
F_N^a(s,t)=F_N^a(s,s) + \oneN \sum_{i=1}^N \E [ G^i(\Bx_s) 
\E (B^i_t-B^i_s | \Gga_s)] = F_N^a(s,s) 
$$
(recall the independence of $\BB_t-\BB_s$ and $\Gga_s$). By the same 
reasoning also $\chi_N^a(s,0)=0$, 
$E_N^a(s,0)=0$, $E_N^a(s,t)=E_N^a(s,s)$ and $\chi_N^a(s,t)=\chi_N^a(s,s)$
for all $N$ and any $t \geq s \geq 0$ (cf. \req{eq:chidef} and 
\req{eq:dedef}). Likewise, by definition $C_N^a(t,s)=C_N^a(s,t)$.

Let $\Phi(C,\chi,D,E): \Ca ([0,T]^2)^4 \to \Ca ([0,T]^2)^4$ denote the 
difference between the two sides of equations \req{eqC1}--\req{eqE}.
In view of the above boundary and symmetry conditions, 
comparing \req{eqCa}--\req{eq:apxf} with \req{eqC1}--\req{eqE} we
see that $\Phi(C_N^a,\chi_N^a,D_N^a,E_N^a) \to 0$ as $N \to \infty$,
uniformly on $[0,T]^2$. It is not hard to check that $\Phi(\cdot)$ is 
continuous with respect to the topology of uniform convergence, 
hence $\Phi(\cdot)=0$ at 
any limit point of $(C_N^a,\chi_N^a,D_N^a,E_N^a)$, which 
also necessarily satisfies the same boundary and symmetry conditions,
thus completing the proof of Proposition \ref{prop-macro}.
\hfill\qed

As already noted, the first step in proving Proposition \ref{comp1} is,
\begin{lem}\label{kval}
Let $\E_{\BJ}$ denotes the expectation with respect to the 
Gaussian law $\P_{\BJ}$ of the disorder $\BJ$. Then, 
for each continuous path $\Bx \in \Ca(\R_+,\R^N)$ and all
$s,t \in [0,T]$ and $i,j \in \{1,\ldots,N\}$,
\begin{equation}\label{eq:kval}
k^{ij}_{ts}(\Bx) :=
\E_{\BJ} [G^i(\Bx_t)G^j(\Bx_s)]= \frac{x_t^j x_s^i}{N}
\nu''( C_N (s,t) ) + {\bf 1}_{i=j} \nu' (C_N(s,t)) \;.
\end{equation}
\end{lem}  
\prf 
Observe that $G^i(\Bx)$ of \req{eq:gdef} are, for any given $\Bx$,
centered, jointly Gaussian random variables.
Further, 
by our choice of $c(\{i_1,\ldots,i_p\})$ it is not hard to verify that
for any $i,i_1,\ldots,i_{p-1}=1,\ldots,N$
\begin{equation}\label{eq:comb1}
c(\{i,i_1,\ldots,i_{p-1}\}) 
 \sum_{1\le j_1,\ldots, j_{p-1} \le N} 
{\bf 1}_{\{i,i_1,\ldots,i_{p-1}\}= \{j,j_1,\ldots,j_{p-1}\}} = 
(p-1)! ({\bf 1}_{j=i} + \sum_{r=1}^{p-1} {\bf 1}_{j=i_r}) \,.
\end{equation}
Hence, by \req{eq:gdef} and \req{eq:comb1}, we have that 
for any given vectors $\Bx$ and $\By$,
\begin{eqnarray*}
\E[G^i(\Bx)G^i(\By)]
&=& \sum_{ p=1}^m \left(\frac{a_p}{(p-1)!}\right)^2 
 \sum_{i_1,\cdots,i_{p-1},j_1,\cdots,j_{p-1}}
\E( J_{i i_1 \cdots i_{p-1}}J_{i j_1 \cdots j_{p-1}} ) 
x^{i_1}\ldots  x^{i_{p-1}}y^{j_1}\ldots  y^{j_{p-1}}\\
&=& \sum_{p=1}^m \frac{a_p^2}{(p-1)!} N^{-(p-1)} \sum_{i_1,\cdots,i_{p-1}}
(1 + \sum_{r=1}^{p-1} {\bf 1}_{i=i_r})
x^{i_1}\ldots  x^{i_{p-1}}y^{i_1}\ldots  y^{i_{p-1}}\\
&=&\sum_{p=1}^m \frac{a_p^2}{(p-1)!}\Big(\oneN 
\sum_{\ell=1}^N x^\ell y^\ell\Big)^{p-1} +
\frac{x^i y^i}{N} \sum_{p=2}^m \frac{a_p^2}{(p-2)!}\Big(\oneN
\sum_{\ell=1}^N x^\ell y^\ell\Big)^{p-2}\,,
\end{eqnarray*}
so that
with $\nu(\rho)=\sum_{p=1}^m \frac{a_p^2}{p!} \rho^p$ we have,
\begin{equation}\label{eq:diagcov}
\E[G^i(\Bx)G^i(\By)]=\nu'\Big( \oneN \sum_{\ell=1}^N
x^\ell y^\ell\Big)
+ \frac{x^i y^i}{N}
\nu''\Big( \oneN \sum_{\ell=1}^N x^\ell y^\ell\Big)\,.
\end{equation}
Further, if $i\neq j$, then by \req{eq:comb1} we also have that
\begin{eqnarray*}
x^i y^j \E[G^i(\Bx)G^j(\By)]&=&
\sum_{p=2}^m \left(\frac{a_p}{(p-1)!}\right)^2
 \sum_{i_1,\cdots,i_{p-1},j_1,\cdots,j_{p-1}}
\E( J_{i i_1 \cdots i_{p-1}}J_{j j_1 \cdots j_{p-1}} ) 
x^i x^{i_1}\ldots  x^{i_{p-1}} y^j y^{j_1}\ldots  y^{j_{p-1}}\\
&=& 
\sum_{p=2}^m \frac{a_p^2}{(p-1)!} N^{-(p-1)} \sum_{i_1,\cdots,i_{p-1}}
\sum_{r=1}^{p-1} {\bf 1}_{j=i_r}
x^i x^{i_1}\ldots  x^{i_{p-1}} y^i y^{i_1}\ldots  y^{i_{p-1}}\\
&=& \frac{x^i y^i x^j y^j}{N}\sum_{p=2}^m \frac{a_p^2}{(p-2)!}
\Big(\oneN \sum_{\ell=1}^N x^\ell y^\ell\Big)^{p-2}\,,
\end{eqnarray*}
implying that when $i \neq j$, 
\begin{equation}\label{eq:offdcov}
\E[G^i(\Bx)G^j(\By)]=\frac{x^j y^i}{N}
\nu''\Big( \oneN \sum_{\ell=1}^N x^\ell y^\ell\Big)\,,
\end{equation}
so replacing $\Bx$ and $\By$ by $\Bx_t$ and $\Bx_s$ respectively,
we immediately get \req{eq:kval} out of \req{eq:diagcov},
\req{eq:offdcov} and the definition of $C_N(\cdot,\cdot)$.
\qed

\medskip
\noindent
{\bf Proof of Proposition \ref{comp1}:}
Fixing a continuous path $\Bx$, 
let $k_t$ denote the operator 
on $L_2(\{1,\cdots N\}\ts[0,t])$
with the kernel $k=k(\Bx)$ of \req{eq:kval}. 
That is, for $f\in L_2(\{1,\cdots N\}\ts[0,t])$,
$u\le t$, $i\in\{1,\cdots,N\}$
\begin{equation}\label{eq:kfdef}
[k_t f]_u^i =\sum_{j=1}^N \int_0^t k_{uv}^{ij} f^j_v dv,
\end{equation}
which is clearly also in $L_2(\{1,\cdots N\}\ts[0,t])$.
We next extend the definition \req{eq:kfdef} to the
stochastic integrals of the form 
$$
[k_t \circ dZ]_u^i=\sum_{j=1}^N \int_0^t k_{uv}^{ij} dZ^j_v,
$$
where $Z^j_v$ is a continuous semi-martingale with respect to 
the filtration $\Fa_t=\sigma(\Bx_u : 0 \leq u \leq t )$ and is
composed for each $j$, of a squared-integrable continuous martingale 
and a continuous, adapted, squared-integrable finite variation
part. In doing so, recall that by \req{eq:kval}
each $k_{uv}^{ij}(\Bx)$ is the finite sum of terms such as
$x^{i_1}_{u}\cdots x^{i_a}_{u}  x^{j_1}_v\cdots x^{j_b}_v$,
where in each term $a$, $b$ and $i_1,\ldots,i_a,j_1,\ldots,j_b$ are
some non-random integers. Keeping for simplicity 
the implicit notation $\int_0^t k_{uv}^{ij} dZ^j_v$ we thus adopt
hereafter the convention of accordingly decomposing such integral
to a finite sum, taking for 
each of its terms the variable $x^{i_1}_{u}\cdots x^{i_a}_{u}$
outside the integral, resulting with the usual It\^o  
adapted stochastic integrals. The latter are well defined, with
$[k_t \circ dZ]_u^i$ being in $L_2(\{1,\cdots N\}\ts[0,t])$ 
(recall Proposition \ref{existN} that $\Bx_t$ has uniformly bounded 
finite moments of all orders
under the joint law $\P_\BJ\otimes \P^N_{\bx_0,\BJ}$,
hence so does the kernel $k_{ts}^{ij}(\Bx)$). 

Equipped with these definitions, we next claim that
\begin{lem}\label{condexp}
Fixing $\tau \in \R_+$, let $V_s^i (\Bx)=\E[G^i(\Bx_s)|\Fa_\tau]$ 
and $Z_s^i (\Bx)=\E[B^i_s|\Fa_\tau]$ for $s \in [0,\tau]$. Then, 
under $\P_\BJ\otimes \P^N_{\bx_0,\BJ}$ we can choose a version of 
these conditional expectations such that the stochastic processes
\begin{eqnarray}
U^i_s(\Bx)&=& x^i_s-x^i_0+\int_0^s f'(K_N(u)) x^i_u du \label{eq:Udef} \\
Z^i_s(\Bx)&=& U^i_s(\Bx) - \int_0^s V^i_u (\Bx) du \;,
\label{eq:UBdef}
\end{eqnarray}
are both continuous semi-martingales with respect to 
the filtration $\Fa_s$, composed of squared-integrable continuous 
martingales and finite variation parts. 
Moreover, such choice satisfies for any $i$ and $s \in [0,\tau]$,
\begin{equation}\label{eq:Vid}
V^i_s + [k_\tau V]^i_s = [k_\tau \circ dU]^i_s\,,
\end{equation}
and $V^i_s = [k_\tau \circ dZ]_s^i$ for any $i$ and all $s \leq \tau$.
Further, for any $u,v \in [0,\tau]$ and $i,j \leq N$, let
\begin{equation}\label{eq:gammadef}
\Gamma_{uv}^{ij} (\Bx) := \E \Big[(G^i(\Bx_u) - V_u^i(\Bx))(G^j (\Bx_v)- 
V_v^j(\Bx))|\Fa_\tau\Big]
\end{equation}
Then, we can choose a version of $\Gamma_{uv}^{il}$ 
such that for any $s,v \leq \tau$ and all $i,l \leq N$, 
\begin{equation}\label{eq:gaid}
\sum_{j=1}^N \int_0^\tau k_{su}^{ij} \Gamma_{uv}^{jl} du 
+ \Gamma_{sv}^{il} = k_{sv}^{il} \;.
\end{equation}
\end{lem}
\prf Since $U^i_s(\Bx)=\int_0^s G^i(\Bx_u) du + B_s^i$ (see \req{eq:integ}),
the relation \req{eq:UBdef} between $\E[B^i_s|\Fa_\tau]$, $U^i_s(\Bx)$ and 
$\E[G^i(\Bx_u)|\Fa_\tau]$ follows, as well as the continuity 
and integrability properties of the semi-martingales $U_s$ and $Z_s$.
Using hereafter $G^i_s$ to denote $G^i(\Bx_s)$, let
\begin{equation}\label{eq:lndf}
\L^N_\tau=\exp\Big\{ \sum_{i=1}^N \int_0^\tau G^i_s dU^i_s(\Bx)
-\half\sum_{i=1}^N \int_0^\tau (G^i_s)^2 ds\Big\} \,.
\end{equation}
By Girsanov formula we have that the  restriction to $\Fa_\tau$ satisfies,
$$\P^N_{\Bx_0,\BJ}|_{\Fa_\tau}=\L^N_\tau \P^N_{\Bx_0,0}|_{\Fa_\tau}$$
Hence, with $\tau \ge s$,
for any bounded $\Fa_\tau$-measurable random variable $\Phi$,
\begin{equation}\label{eq:conde}
\E [ G^i_s \Phi ] = \E_\BJ\E_{\P^N_{\Bx_0,\BJ}}[G^i_s \Phi]=
\E_{\P^N_{\Bx_0,0}}[\E_\BJ[G^i_s \L^N_\tau]\Phi]
=\E \Big[\frac{\E_\BJ[G^i_s \L^N_\tau]}{\E_\BJ[\L^N_\tau]}\Phi\Big]\,,
\end{equation}
where the right-most identity is due to the change of measure formula 
$\Q^N_{\bx_0} = \E_\BJ(\L^N_\tau) \P^N_{\bx_0,0}$ for the
annealed law $\Q^N_{\bx_0}=\E_\BJ \P^N_{\bx_0,\BJ}$, restricted to
$\Fa_\tau$.
With \req{eq:conde} holding for all bounded $\Fa_\tau$-measurable $\Phi$,
it follows that,
$$
V_s^i = \E[G^i_s|\Fa_\tau]
=\frac{\E_\BJ[G^i_s \L^N_\tau]}{\E_\BJ[\L^N_\tau]}\,,
$$
and the identity \req{eq:Vid} follows by the 
Gaussian change of measure identity \req{eq:mean}
of Proposition \ref{prop-alice}. Exactly the same line of reasoning
shows that,
$$
\Gamma_{uv}^{ij} = 
\frac{\E_\BJ[(G^i_u-V^i_u)(G^j_v-V^j_v) \L^N_\tau]}{\E_\BJ[\L^N_\tau]}\,,
$$
and the identity \req{eq:gaid} follows by the 
identity \req{eq:cov} of Proposition \ref{prop-alice}.
Noting that $dZ=dU-V$ (see \req{eq:UBdef}), we have by \req{eq:Vid}
that for all $i$ and $s \in [0,\tau]$,
$$
[k_\tau \circ dZ]_s^i = [k_\tau \circ dU]_s^i - [k_\tau V]_s^i = V_s^i
$$
as claimed.
\hfill\qed

We now apply \req{eq:Vid} to derive \req{eq:apxa},
the easy part of Proposition \ref{comp1}. To this end,
fixing $s,t \in [0,T]^2$ let 
$$
\widehat{A}_N(t,s)=\oneN \sum_{i=1}^N V^i_t(\Bx) x^i_s \,,
$$
for $\tau=t \vee s$, noting that since $x^i_s$ is measurable on
$\Fa_\tau$, by \req{eq:hfdef},
$$
A_N^a (t,s) =\E \Big[ \oneN \sum_{i=1}^N \E [G^i_t x^i_s|\Fa_\tau] \Big]
= \E [ \widehat{A}_N(t,s) ] = \widehat{A}_N^a (t,s) \;.
$$
With $t \leq \tau$, combining \req{eq:Vid} and \req{eq:Udef} we 
get 
 that
$$
\widehat{A}_N(t,s) + \oneN \sum_{i,j=1}^N \int_0^\tau 
x_s^i k_{tu}^{ij} V_u^j du = 
\oneN \sum_{i,j=1}^N \int_0^\tau f'(K_N(u)) x_s^i k_{tu}^{ij} x_u^j du 
+ \oneN \sum_{i,j=1}^N \int_0^\tau x_s^i k_{tu}^{ij} dx_u^j 
$$
(suppressing the dependence of $k^{ij}_{tu}$ and
$V_u^j$ on $\Bx$, and following our convention regarding 
stochastic integrals such as $\int_0^\tau x_s^i k_{tu}^{ij} dx_u^j$).
Using the explicit expression of 
$k^{ij}_{tu}(\Bx)$, and collecting terms while       
changing the order of summation and integration, we arrive
at the identity,
\begin{eqnarray}\label{eq:Ahatid}
&& \widehat{A}_N(t,s) + \int_0^\tau 
C_N(s,u) \nu''(C_N(t,u)) \widehat{A}_N (u,t) du 
+ \int_0^\tau \nu'(C_N(t,u)) \widehat{A}_N (u,s) du 
\nonumber \\ 
&=& 
\int_0^\tau f'(K_N(u)) C_N(s,u) \nu''(C_N(t,u)) C_N (u,t) du 
+ \int_0^\tau f'(K_N(u)) \nu'(C_N(t,u)) C_N (u,s) du 
\nonumber \\
&+& \int_0^\tau C_N(s,u) \nu''(C_N(t,u)) d_u C_N(u,t) 
+ \int_0^\tau \nu'(C_N(t,u)) d_u C_N (u,s) \,.
\end{eqnarray}
Applying Lemma \ref{ito} for the semi-martingales $x=y=z=w=\Bx$, 
and the polynomials $Q(\rho)=\nu'(\rho)$ and $P(\rho)=\rho$ evaluated
at $\sigma=\tau$, $\theta=s$ and $v=t$, we replace the
stochastic integrals of \req{eq:Ahatid} by
\begin{eqnarray}\label{eq:nostoch}
&& \nu'(C_N(\tau,t)) C_N(\tau,s) - \nu'(C_N(0,t)) C_N(0,s) \nonumber \\
&-& \frac{1}{2N} C_N(t,t) \int_0^\tau C_N(u,s) \nu''(C_N(u,t)) du 
-\frac{1}{N} C_N(s,t) \int_0^\tau \nu'(C_N(u,t)) du \;.
\end{eqnarray}
Clearly, $\widehat{A}_N(t,s) = \E [ A_N(t,s) | \Fa_\tau ]$ 
has the same uniform moment bounds 
of Proposition \ref{tightchic} as $A_N$ and further 
inherits the self-averaging property \req{eq:l2self} from $A_N$.
Hence, 
we may and shall apply Corollary \ref{cor-self}
with possibly $\widehat{A}_N$ as one of the arguments of the
locally Lipschitz function $\Psi(z)$ of at most polynomial growth 
at infinity. Doing so for the functions $z_1 z_2 \nu''(z_3)$,
and $z_1 \nu'(z_2)$ and applying Corollary \ref{cor-self}
also for $f'(z_1) z_2 \nu''(z_3) z_3$ and
$f'(z_1) \nu'(z_2) z_3$, upon utilizing the uniform convergence
with respect to the points $(s_j,t_j) \in [0,T]^2$, we deduce
from \req{eq:Ahatid} and \req{eq:nostoch} that
\begin{eqnarray*}
\widehat{A}^a_N(t,s) &+& \int_0^\tau 
 C_N^a(s,u) \nu''(C_N^a(t,u)) \widehat{A}_N^a (u,t) du 
+ \int_0^\tau \nu'(C_N^a (t,u)) \widehat{A}_N^a (u,s) du 
\nonumber \\ 
&\simeq&
\int_0^\tau f'(K_N^a(u)) C_N^a(s,u) \nu''(C_N^a(t,u)) C_N^a (u,t) du 
\nonumber \\ 
&+& 
\int_0^\tau f'(K_N^a(u)) \nu'(C_N^a(t,u)) C_N^a (u,s) du 
+ \nu'(C_N^a(\tau,t)) C_N^a(\tau,s) - \nu'(C_N^a(0,t)) C_N^a(0,s) \,.
\end{eqnarray*}
Finally, recall that 
\begin{equation}\label{sa2}
\widehat{A}_N^a(t,s)=A_N^a(t,s)=D_N^a(s,t) + f'(K_N^a(t)) C_N^a(s,t)\,, 
\end{equation}
with the corresponding replacement for $\widehat{A}_N^a(u,t)$ and 
$\widehat{A}_N^a(u,s)$. With $\tau=t \vee s$, 
we indeed arrive at \req{eq:apxa}.

We now turn to the more involved part of 
Proposition \ref{comp1}, namely, the derivation of \req{eq:apxf}.
To this end, as we have seen already, it suffices to consider
$s \geq t$, as we do hereafter. To this end, taking $\tau=s$ 
we use the notation
$V_u=\E[G_u|\Fa_s]$ of Lemma \ref{condexp} while suppressing the 
dependence on $\Bx$. Since $\BB_t=U_t - \int_0^t G_v dv$ with
$U_t$ being $\Fa_s$-measurable, we deduce that  
$$ 
\E[G^i_s B^i_t ] = \E[ G^i_s U^i_t] - \int_0^t \E [G^i_v G^i_s] dv
=\E \Big[ \E[ G^i_s |\Fa_s](U^i_t-\int_0^t V^i_v dv) 
-\int_0^t \Gamma_{sv}^{ii} dv \Big] \,,
$$
(recall that $\E(G_v-V_v|\Fa_s)=0$ for all $v \leq s$
and $\Gamma_{uv}^{ij} := \E [(G_u^i - V_u^i)(G_v^j - V_v^j)|\Fa_s]$
is per \req{eq:gammadef}). Further, by
\req{eq:UBdef} and the identity 
$\E[ G^i_s |\Fa_s]=[k_s \circ dZ]_s^i$ of Lemma \ref{condexp},
we get
\begin{equation}\label{eq2}
\E[G^i_s B^i_t ] +\E [ \int_0^t \Gamma_{sv}^{ii} dv ]
= \E ( [k_s \circ dZ]_s^i Z^i_t ) \,.
\end{equation}
Since $Z^i_t=\E[B^i_t|\Fa_s]$ and 
$[k_s \circ dZ]_s^i$ is $\Fa_s$-measurable, we have that
$$
\E ( [k_s \circ dZ]_s^i Z^i_t ) = \E ( [k_s \circ dZ]_s^i B^i_t ) 
=  \E ( \E( [k_s \circ d\BB]_s^i |\Fa_s ) B^i_t ) \,,
$$
where the right-most identity holds since the kernel of the 
linear operator $k_s$ is $\Fa_s$-measurable and 
$\E[\BB_u|\Fa_s]=Z_u$ for all $u \leq s$ (recall that
$[k_s \circ dZ]_s^i$ is the $L_2$-limit of 
discrete sums with mash size going to zero).
In view of \req{eq:integ} 
and \req{eq:kval} we have that
\begin{eqnarray}\label{eq:kdb}
[k_s \circ d\BB]_s^i &=& [k_s \circ d\Bx]_s^i + [k_s f'(K_N) \Bx]_s^i 
- [k_s G]_s^i \nonumber \\
&=&  \int_0^s \nu'( C_N (s,u)) dx_u^i 
+ \int_0^s \nu''( C_N(s,u)) x_u^i d_u C_N(s,u) \\
&+& \int_0^s \nu'( C_N(s,u)) f'(K_N(u)) x_u^i du 
+ \int_0^s \nu''( C_N(s,u)) f'(K_N(u)) C_N(s,u) x_u^i du \nonumber \\
&-& \int_0^s \nu''( C_N(s,u)) x_u^i A_N(u,s) du 
- \int_0^s \nu'( C_N(s,u)) G_u^i du 
\nonumber
\end{eqnarray}
Using It\^o's formula for $x_u^i \nu'(C_N(s,u))$ we replace
the two stochastic integrals in \req{eq:kdb} by
\begin{equation}\label{nostoch2}
x_s^i \nu'(C_N(s,s)) - x_0^i \nu'(C_N(s,u)) -
\frac{1}{2N} C_N(s,s) \int_0^s \nu'''(C_N(s,u)) x^i_u du 
-\oneN x^i_s \int_0^s \nu''(C_N(s,u)) du \;.
\end{equation}
Recall that by \req{eq:hfdef} and \req{eq2}
$$
F_N^a (s,t) + \E [ \oneN \sum_{i=1}^N \int_0^t \Gamma_{sv}^{ii} dv ] =
\E ( \oneN \sum_{i=1}^N \E( [k_s \circ d\BB]_s^i |\Fa_s ) B^i_t ) \,,
$$
which by the preceding is the expectation of 
\begin{eqnarray}\label{eq:detail}
&& \chi_N(s,t) \nu'(C_N(s,s)) - \chi_N(0,t) \nu'(C_N(s,0)) \nonumber \\
&-& \frac{1}{2N} C_N(s,s) \int_0^s \nu'''(C_N(s,u)) \chi_N(u,t) du 
-\oneN \int_0^s \nu''(C_N(s,u)) \chi_N(s,t) du \nonumber \\ 
&+& \int_0^s f'(K_N(u)) \nu'( C_N(s,u)) \chi_N(u,t) du 
+ \int_0^s f'(K_N(u)) \nu''( C_N(s,u)) C_N(s,u) \chi_N(u,t) du \nonumber \\
&-& \int_0^s \nu''( C_N(s,u)) \chi_N(u,t) \widehat{A}_N (u,s) du 
- \int_0^s \nu'( C_N(s,u)) F_N(u,t) du +\kappa_N(s,t) \,,
\end{eqnarray}
where in view of \req{eq:hfdef}
$$
\kappa_N(s,t) := \oneN \sum_{i=1}^N \int_0^s \nu'( C_N(s,u)) 
(G_u^i-V_u^i) B^i_t du \,.
$$
Recall \req{eq:UBdef} that
$B_t^i=Z_t^i-\int_0^t (G_v^i-V_v^i) dv$, hence 
$$
\kappa_N(s,t) = \oneN \sum_{i=1}^N \int_0^s \nu'( C_N(s,u)) 
(G_u^i-V_u^i) Z^i_t du 
- \oneN \sum_{i=1}^N 
\int_0^s \nu'( C_N(s,u)) \int_0^t (G_u^i-V_u^i) (G_v^i-V_v^i) dv du \;.
$$
As both $\nu'(C_N(u,s))$ and $Z_i^t$ are $\Fa_s$-measurable while
$\E(G_u-V_u|\Fa_s)=0$ for all $u \leq s$, the expectation of the
first term on the right-side vanishes. Further, conditioning 
on $\Fa_s$, we have in view of \req{eq:gammadef} that
\begin{equation}\label{eq:kaexp}
\E[\kappa_N(s,t)] = -\E\Big[ \oneN \sum_{i=1}^N \int_0^s \nu'( C_N(s,u)) 
\int_0^t \Gamma_{uv}^{ii} dv du \Big] \;.
\end{equation}
All terms of \req{eq:detail}
apart from $\kappa_N$ are of the form covered by 
Corollary \ref{cor-self} for functions $\Psi(z)$ similar to 
those encountered in the derivation of \req{eq:apxa},
(namely, $z_1 \nu'(z_2)$, $z_1 \nu'''(z_2) z_3$, $z_1 \nu''(z_2)$,
$f'(z_1) \nu'(z_2) z_3$, $f'(z_1) z_2 \nu''(z_3) z_3$ 
and $z_1 z_2 \nu''(z_3)$). Utilizing their uniform convergence and
\req{eq:kaexp}, while recalling that $\chi_N^a(0,t)=0$,
\req{sa2} and the analogous 
$$
F_N^a(u,t)=E_N^a(u,t) + f'(K_N^a(u)) C_N^a(s,u)\,, 
$$
is not hard to check that we get \req{eq:apxf}, once we prove the
following lemma.
\begin{lem}\label{cancelem}
For $v \in [0,s]$, let 
$$
\phi_N(s,v) := \nu'(C_N(s,v))  - \oneN \sum_{i=1}^N \Gamma_{sv}^{ii} 
- \oneN \sum_{i=1}^N \int_0^s \nu'( C_N(s,u)) \Gamma_{uv}^{ii} du 
$$
(where $\Gamma$ is defined with $\tau=s$). Then, for any $t \leq s$,
$$
\int_0^t \E[ \phi_N(s,v) ]  dv \simeq 0
$$
\end{lem}
\prf Fixing $v \leq s$, recall that by \req{eq:kval} we have that 
$$
\phi_N(s,v) +\oneN \nu''(C_N(s,v)) C_N(s,v)
- \frac{1}{N^2} \sum_{i,j=1}^N \int_0^s \nu''(C_N(s,u)) x^j_s x^i_u
\Gamma_{uv}^{ji} du = \oneN \sum_{i=1}^N \Big[ k_{sv}^{ii} - \Gamma_{sv}^{ii} 
- \sum_{j=1}^N \int_0^s k_{su}^{ij} \Gamma_{uv}^{ji} du \Big]
$$
Further, since we had $\tau=s$, the right-side is identically zero by
\req{eq:gaid}, implying by the $\Fa_s$-measurability of 
$\Bx_s$, $\Bx_u$ and some algebraic manipulations that
\begin{eqnarray*}
&& \phi_N(s,v) +\oneN \nu''(C_N(s,v)) C_N(s,v) = 
\int_0^s \nu''(C_N(s,u)) \frac{1}{N^2} \sum_{i,j=1}^N 
x^j_s x^i_u \E[ (G^j_u-V^j_u)(G^i_v-V^i_v) |\Fa_s ] du \\
&=& 
\int_0^s \nu''(C_N(s,u)) 
\E\Big[ (A_N(u,s) - \widehat{A}_N (u,s)) (A_N(v,u) - \widehat{A}_N (v,u))
 |\Fa_s \Big] du \;. 
\end{eqnarray*}
Consequently, with $C_N(s,u)$ measurable on $\Fa_s$ we have that
$$
\E [ \phi_N(s,v) ] = 
\int_0^s \E \big[ \nu''(C_N(s,u)) 
 (A_N(u,s) - \widehat{A}_N (u,s)) (A_N(v,u) - \widehat{A}_N (v,u))
 \big] du -\oneN \E [ \nu''(C_N(s,v)) C_N(s,v) ]  
\;. 
$$
In view of Proposition \ref{tightchic} and \req{eq:l2self},
the first term converges to zero uniformly in $s,v$ by the 
uniform $L_2$ convergence of $A_N-A_N^a$ (hence also of
$A_N-\widehat{A}_N$), and the uniform (in $[0,T]^2$ and
$N$) bound on each moment of $C_N$ and $A_N$ (hence on
those of $\widehat{A}_N$ as well). By same reasoning,
the second term converges to zero at rate $1/N$, uniformly
in $s,v$. Utilizing the uniformity of the convergence, we
see that 
$\int_0^t \E[ \phi_N(s,v) ]  dv \simeq 0$ as claimed.
\hfill\qed

\section{Differentiability and uniqueness for 
the limiting dynamics}\label{sec-uniq}

We start the proof of Theorem \ref{thm-macro}
by the next lemma 
relating the solutions of \req{eqC1}--\req{eqE} with those
of \req{eqR}--\req{eqZ}.
\begin{lem}\label{lem-diff}
Fixing $T<\infty$, suppose $(C,\chi,D,E)$ is a solution 
of the integral equations \req{eqC1}--\req{eqE}
in the space of bounded continuous functions on $[0,T]^2$
subject to the symmetry condition $C(s,t)=C(t,s)$ and the
boundary conditions $E(s,0)=0$ for all $s$, and
$E(s,t)=E(s,s)$ for all $t \geq s$. 
Then, $\chi(s,t)=\int_0^t R(s,u)du$ where
$R(s,t)=0$ for $t>s$, $R(s,s)=1$ and
for $T \geq s>t$, the bounded and
absolutely continuous functions 
$C$, $R$ and $K(s)=C(s,s)$ necessarily satisfy 
the integro-differential equations \req{eqR}--\req{eqZ}.
\end{lem}
\prf Consider the integral operator $k_C$ 
on $\Ca ([0,T])$ given by,
$$
[k_{C} h](s):=-\int_0^s \nu'(C(s,u)) h(u) du \,,
$$
and let
$$
h(s,t):=-f'(C(s,s))\chi(s,t)
-\int_0^s\chi(u,t)\nu''(C(s,u))D(s,u) du+
\chi(s,t) \nu'(C(s,s)) -\int_0^{t \wedge s} \nu'(C(s,u)) du
$$
Then, per fixed $t$, the equation \req{eqE} states
that $E(s,t)=[k_C E(\cdot,t)](s) + h(s,t)$. Since the (continuous)
kernel $\nu'(C(s,u))$ of $k_C$ is uniformly bounded on $[0,T]^2$, 
it follows by Picard iterations (splitting $[0,T]$ to sufficiently
small time intervals to guarantee convergence of the series
$\sum_n k_C^n$), that
\begin{equation}
E(s,t)
=\sum_{n\ge 0} [k_{C}^n h(.,t)](s)
= h(s,t) + \int_0^s \kappa_C (s,v) h(v,t) dv \,,
\label{blurp}
\end{equation}
with a uniformly bounded kernel $\kappa_{C}$.
Plugging \req{blurp} into \req{eq:chi}, we find by Fubini's theorem 
that
$$
\chi(s,t)= 
s \wedge t + \int_0^s[\int_0^{t \wedge v} \nu'(C(v,u)) du]\kappa_1(s,v) dv
+ \int_0^s \chi(v,t)\kappa_2(s,v) dv \,,
$$
for some uniformly bounded functions $\kappa_1$ and
$\kappa_2$ which depend only on $C$ and $D$.
Applying Picard's iterations once more, now with respect to the
integral operator $[\kappa_2 g](s) = \int_0^s \kappa_2(s,v) g(v) dv$,
we deduce that for some uniformly bounded $\kappa_3$ and $\kappa_4$,
$$
\chi(s,t)= s\wedge t + \int_0^s \Big[ (u\wedge t) \kappa_3 (s,u) +
\big[ \int_0^{t \wedge u} \nu'(C(u,v)) dv \big] \kappa_4(s,u) \Big] du \;.
$$
On $s\ge t,$ the function $s\wedge t=t$ is continuously 
differentiable, hence we easily conclude by Fubini's theorem 
that $t\ra \chi(s,t)$ is continuously differentiable on $s\ge t$,
with $\chi(s,t)=\int_0^t R(s,u) du$ for the bounded
continuous function
$$
R(s,t)=1+\int_t^s [\kappa_3 (s,u) + \nu'(C(u,t)) \kappa_4(s,u)]du \;.
$$
In particular, $R(s,s)=1$ for all $s$.
The condition $E(s,t)=E(s,s)$ for $s \geq t$ implies
by \req{eq:chi} that $\chi(s,t)=\chi(s,s)$ for $s \geq t$.
Similarly, with $E(s,0)=0$, it follows that $\chi(s,0)=0$ 
for all $s$. From \req{eqC1} we have that $C(s,t)-\chi(s,t)$
is differentiable with respect to its 
second argument $t$, with a bounded, continuous derivative
$D=\partial_2 (C-\chi)$. Consequently, $\partial_2 C = D + R$
where $R(s,t) = (\partial_2 \chi) (s,t) = 0$ for all $t>s$ 
due to the boundary condition $\chi(s,t)=\chi(s,s)$. 
Further, $C(s,t)=C(t,s)$ implying that
$\partial_1 C (s,t) = \partial_2 C(t,s) = D(t,s)+R(t,s)$
on $[0,T]^2$. Thus, combining the identity
\begin{eqnarray*}
C(s,t \vee s)\nu'(C(t \vee s,t))-C(s,0)\nu'(C(0,t))&=&
\int_0^{t \vee s} \nu'(C(t,u)) (\partial_2 C) (s,u) du \\ 
&+&\int_0^{t \vee s} C(s,u) \nu''(C(t,u)) (\partial_2 C) (t,u) du \,,
\end{eqnarray*}
with \req{eqD} we have that for all $t,s \in [0,T]^2$,
\begin{equation}\label{eq:tempD}
D(s,t)=-f'(K(t)) C(t,s)+
\int_0^{t \vee s} \nu'(C(t,u))R(s,u) du
+\int_0^{t \vee s} C(s,u)\nu''(C(t,u)) R(t,u) du \,.
\end{equation}
Interchanging $t$ and $s$ in \req{eq:tempD} and adding
$R(t,s)=0$ when $s>t$, results for $s>t$ with
$$
(\partial_1 C) (s,t) =
-f'(K(s)) C(s,t)+
\int_0^s\nu'(C(s,u))R(t,u) du
+\int_0^s C(t,u)\nu''(C(s,u)) R(s,u) du \,,
$$
which is \req{eqC} for $\beta=1$. 

Since $K(s)=C(s,s)$, with $C(s,t)=C(t,s)$ 
and $\partial_2 C = D+R$, it follows that for all $h>0$,
$$
K(s)-K(s-h) = \int_{s-h}^{s}(D(s,u) + R(s,u))du
+ \int_{s-h}^{s}(D(s-h,u) + R(s-h,u))du \,.
$$
Recall that $R(s,u)=0$ for $u > s$, hence, 
dividing by $h$ and taking $h \downarrow 0$, we thus get
by the continuity of $D$ and that of 
$R$ for $s \geq t$ that $K(\cdot)$ is differentiable,
with $\partial_s K(s)=2 D(s,s)+R(s,s)=2 D(s,s) + 1$,
resulting by \req{eq:tempD} with \req{eqZ} for $\beta=1$.

Further, it follows from \req{eq:chi} that
$(\partial_1 \chi) (u,t) = E(u,t)+ 1_{u<t}$.
Hence, combining the identity
$$ 
\chi(s,t)\nu'(C(s,s))-\chi(0,t)\nu'(C(s,0)) =
\int_0^s \nu'(C(s,u)) (\partial_1 \chi) (u,t) du 
+\int_0^s \chi(u,t) \nu''(C(s,u)) (\partial_2 C) (s,u) du \,,
$$
with \req{eqE} we have that for all $T \geq s \geq t$,
\begin{equation}\label{eq:tempE}
E(s,t)=-f'(K(s)) \chi (s,t)+
\int_0^s \chi (u,t)\nu''(C(s,u)) R(s,u) du 
\end{equation}
(recall that $\chi(0,t)=\chi(0,0)=0$). Let
\begin{equation}\label{eq:smE}
g(s,t):=-f'(K(s)) R (s,t)+\int_0^s R(u,t)\nu''(C(s,u)) R(s,u) du \,,
\end{equation}
for $s,t \in [0,T]^2$. Recall that $\chi(s,t)=\int_0^t R(s,v) dv$, 
so by Fubini's theorem, \req{eq:tempE} amounts to 
$E(s,t)=\int_0^t g(s,v) dv$ for all $s \geq t$. Further, with 
$E(s,t)=E(s,s)$ when $t>s$, it follows that 
$$
E(s,t)=\int_0^{t \wedge s} g(s,v) dv
$$ 
for all $s,t \leq T$. Putting this into \req{eq:chi} we have 
by yet another application of Fubini's theorem that 
$$
\int_0^t R(s,u) du = \chi(s,t) 
= t + \int_0^s \int_0^{t \wedge u} g(u,v) dv du
= t + \int_0^t \int_v^s g(u,v) du dv \,,
$$
for any $s \geq t$. Consequently, for every $t \leq s$,   
$$
R(s,t) = 1 + \int_t^s g(u,t) du \,,
$$
implying that $\partial_1 R = g$ for a.e. $s>t$, 
which in view of \req{eq:smE} gives \req{eqR} for $\beta=1$, 
thus completing the proof of the lemma. 
\qed

We proceed by showing that the system of equations of 
interest to us admits at most one solution.
\begin{prop}\label{uniqueness}
Let $T\ge 0$ and $\D_T=\{s,t\in (\R^+)^2 : 0\le t\le s\le T\}$.
There exists at most one solution $(R,C,K)
\in \Ca_b^1( \D_T)\ts \Ca^1_b(\D_T)\ts \Ca_b^1([0,T])$
 to (S):=(\ref{eqR},\ref{eqC},\ref{eqZ}) with $C(s,t)=C(t,s)$
and boundary conditions
\begin{eqnarray}
R(s,s)&\equiv& 1\quad\forall s\ge 0\label{bcR}\\ 
C(s,s)&=&K(s) \quad\forall s\ge 0\label{bcC}\\
C(0,0)&=&K(0)\quad \mbox{ known.}\nonumber
\end{eqnarray}
\end{prop}
\prf 
As mentioned already, we may and shall take 
$\beta=1$ (by scaling $\beta^2 \nu \mapsto \nu$, 
with $\beta^2 \psi \mapsto \psi$ accordingly), just as we 
have done throughout this paper. 

Assume that there are two solutions
$(R,C,K)$ and $(\tilde R,\tilde C,\tilde K)$
of (S) with boundary condition
(BC):=(\ref{bcR},\ref{bcC}). We shall prove by Gronwall's
type argument that these two solutions have to coincide.
To do so we first show that the response function 
$R$ is a Lipschitz function of the covariance
functions $(C,K)$ and then  that the covariances 
obey integro-differential Gronwall type
inequalities. We then use 
Gronwall arguments repeatedly to conclude.
In what follows, $T$ is
fixed and all the constants (which eventually depend
on $T$) will be denoted by $M$, even though
they may change from line to line.

\def\cro{{{\mbox{cr}}}}

\begin{itemize}
\item
{\it $R$ is a Lipschitz function of the covariance}

Let  $(R,C,K)$ be a solution to
(S) and denote 
$$H_C(s,t):=e^{\int_t^s f'(K(u)) du} R(s,t).$$
Then, from \req{eqR} and \req{bcR}, we deduce that
$H$ satisfies
$$\partial_s H_C(s,t)=\int_t^s H_C(u,t)H_C(s,u) \nu''(C(s,u)) 
du,\quad \mbox{for}\quad s\ge t,\qquad H(t,t)=1.$$
Note that $H$ only depends on $C$.
This equation was studied in \cite{GM} where it
was shown that if $\NC_n$ denotes the
set of involutions of $\{1,\cdots,2n\}$
without fixed points and without crossings
and if  $\cro(\sigma)$
is defined to be the set of indices $1\le i\le 2n$
such that $i<\sigma(i)$, $H$ is
given by
\begin{equation}\label{generalformula}
H_C(s,t)=1+\sum_{n\ge 1} 
\int_{t\le t_1\cdots \le t_{2n}\le s}
\sum_{\sigma\in \NC_n}\prod_{i\in
\cro(\s)}  \nu''(C(t_{\s(i)},t_i))dt_1\cdots dt_{2n}
\end{equation}
The number of non-crossing partitions $|\NC_n|$
of $2n$ elements 
is given by the Catalan number $C_n$
which is at most
$2^n$. As 
$\int_{t\le t_1\cdots \le t_{2n}\le s} dt_1\cdots dt_{2n} \leq M/(2n)!$
and $\sup_{(t,s)\in\D_T}\nu''(C(s,t))$ is uniformly bounded
by hypothesis, the above sum converges absolutely. Further,
telescoping each
$\prod_{i\in \cro(\s)}  \nu''(C(t_{\s(i)},t_i))
-\prod_{i\in \cro(\s)}  \nu''(\tilde C(t_{\s(i)},t_i))$, by the
uniform Lipschitz continuity of $\nu''(\cdot)$ on compacts, we thus
find a finite constant $M$ such
that for any pair $C,\tilde C\in \Ca_b(\D_T)$ and any $t,s\in\D_T$,
\begin{equation}\label{ineqH}
|H_C(s,t)-H_{\tilde C}(s,t)|\le
M \int_{t\le t_2\le t_1\le s} |C(t_1,t_2)-\tilde C(t_1,t_2)|
dt_1dt_2.
\end{equation}
Thus, if $(R,C,K)$ and $(\tilde R, \tilde C,\tilde K)$
are two solutions of (S)
in $\Ca_b^1(\D_T)\ts\Ca_b^1(\D_T)\ts\Ca_b^1([0,T])$,
since 
$K$ is uniformly bounded and $f'(\cdot)$ is locally Lipschitz, we obtain
\begin{equation}\label{ineqR}
|R(s,t)-{\tilde R}(s,t)|\le
M \int_{t\le t_2\le t_1\le s} |C(t_1,t_2)-\tilde C(t_1,t_2)|
dt_1dt_2+M\int_t^s |K(u)-\tilde K(u)|du.
\end{equation}
\item {\it Bounds on the difference of the covariances on $s\ge t$}

Integrating (\ref{eqC}) yields for $s\ge t$
\begin{eqnarray*}
C(s,t)&=&K(t)-\int_t^s f'(K(u)) C(u,t) du
+\int_t^sdu\int_0^t dv \nu'(C(u,v)) R(t,v)
 \\
&&+\int_t^sdu\int_0^t dv \nu''(C(u,v))
C(t,v)R(u,v)
+\int_t^sdu\int_t^u dv\nu''(C(u,v))
C(v,t)R(u,v).
\end{eqnarray*}
Hence, if $(R,C,K)$ and $(\tilde R,\tilde C,\tilde K)$
are two solutions of (S),
\begin{eqnarray}
|C-\tilde C|(s,t)&\le &
M \Big[ |K-\tilde K|(t)+\int_t^s |K-\tilde K|(u)du
+\int_t^s|C-\tilde C|(u,t) du
+\int_t^s du\int_0^t dv |C-\tilde C|(u,v)\nonumber\\
&&+\int_t^s du\int_0^t dv|C-\tilde C|(t,v)
+\int_t^s du\int_0^t dv|R-\tilde R|(t,v)
+\int_t^s du\int_0^t dv|R-\tilde R|(u,v)\nonumber\\
&&+\int_t^sdu\int_t^u dv|C-\tilde C|(u,v)+
\int_t^sdu\int_t^u dv|C-\tilde C|(v,t)+
\int_t^sdu\int_t^u dv|R-\tilde R|(u,v)\Big]\nonumber\\
&:=& I_1(s,t)+I_2(s,t)+\cdots +I_{10}(s,t)\label{b1}
\end{eqnarray}

\item{\it Bounds on the differences of the covariances on the diagonal}

Similarly, integrating \req{eqZ} gives
$$K(t)=K(0)-2\int_0^t f'(K(u)) K(u)du +
t+2\int_0^t du \int_0^u dv \psi(C(u,v)) R(u,v) \,, $$
yielding in case $K(0)=\tilde K(0)$ that 
\begin{eqnarray}
|K-\tilde K|(t)&\le &M\Big[\int_0^t |K-\tilde K|(u) du
+\int_0^t du \int_0^u
|C-\tilde C|(u,v) dv +\int_0^t du \int_0^u |R-\tilde R|(u,v) dv \Big]
\label{ineqK}
\end{eqnarray}
Plugging (\ref{ineqR}) into (\ref{ineqK})
yields 
\begin{equation}
|K-\tilde K|(t) \le  M[\int_{0\le t_1\le t_2\le t} |C-\tilde C|(
t_2,t_1) dt_1dt_2+\int_0^t |K-\tilde K|(u) du]\label{ineqK2}
\end{equation}
Recall that by Gronwall's
lemma, if $h,g$ are two non-negative functions
such that
$$h(t)\le g(t)+A\int_0^t h(s)ds$$
for some $A\ge 0$,
then
$$h(t)\le g(t) + A \int_0^t g(v) e^{A(t-v)} dv \le e^{At} g(t)$$
where the last inequality holds
when $g$ is non-decreasing.
Applying this inequality with 
$$g(t)=\int_{0\le t_1\le t_2\le t} |C-\tilde C|(
t_2,t_1) dt_1dt_2
$$ which is non-negative and non-decreasing
yields
\begin{equation}
|K-\tilde K|(t) \le  M[\int_{0\le t_1\le t_2\le t} |C-\tilde C|(
t_2,t_1) dt_1dt_2]\label{ineqK3}
\end{equation}

\item{\it The final fixed point argument}

We now consider
$$D(s):=\int_0^s |C-\tilde C|(s,t) dt,$$
noting that (\ref{ineqR}) and (\ref{ineqK3})
imply that
\begin{equation}\label{ineqR4}
|R(s,t)-{\tilde R}(s,t)|\le
M \int_{0\le t_2\le t_1\le s} |C(t_1,t_2)-\tilde C(t_1,t_2)|
dt_1dt_2 =M \int_{0\le t_1\le s} D(t_1) dt_1 ,
\end{equation}
and 
\begin{equation}\label{ineqK4}
|K-\tilde K|(t) \le  M[\int_{0\le  t_2\le t} D(t_2) dt_2 ].
\end{equation}
Thus, 
 integrating (\ref{b1}) with respect 
to $t$ and 
observing that
\begin{eqnarray*}
\int_0^s ( I_1(s,t)+I_2(s,t))dt
&\le& M\int_0^s D(u) du\mbox{ by (\ref{ineqK4}),}\\
\int_0^s (I_4(s,t)+I_8(s,t))dt &\le& 
M\int_0^s dt \int_t^s du\int_0^u dv|C-\tilde C|(u,v)
\le M \int_0^s D(u) du \mbox{ and }  \\
\int_0^s (I_9(s,t)+I_5(s,t)+I_3(s,t)) dt 
&\le & M \int_0^s D(u) du\quad  \mbox{ by definition of D and Fubini,}\\
\int_0^s (I_6(s,t)+I_7(s,t)+I_{10}(s,t))dt
&\le& M\int_0^s dt\int_t^s du\int_0^u dv|R-\tilde R|(u,v)
\le M\int_0^s D(u)du\mbox{ by (\ref{ineqR4}), }
\end{eqnarray*}
we obtain from \req{b1} that
$$
D(s)\le M \int_0^s D(u) du \,.
$$ 
Recall that $D$ is non-negative and non-decreasing, 
so by the preceding Gronwall argument, now with $g=0$ we
conclude that $D(s)=0$ for all $s\in [0,T]$.
This in turn implies 
by (\ref{ineqR4}) and (\ref{ineqK4})
that
$$K(t)=\tilde K(t),\qquad R(s,t)=\tilde R(s,t)
\quad\mbox{ for all } (t,s)\in\D_T$$
and $C(s,t)=\tilde C(s,t)$ for almost all $t\le s$ and
all $s\le T$. Either by (\ref{b1}) or directly 
by the continuity of the covariances we conclude that
$$ C(s,t)=\tilde C(s,t)
\quad\mbox{ for all } (t,s)\in\D_T$$
which finishes the proof.
\hfill\qed
\end{itemize}

We conclude this section with the,

\medskip
\noindent
{\bf Proof of Theorem \ref{thm-macro}:}
Recall Proposition \ref{tightchic} that we have pre-compactness of 
$(A_N^a,F_N^a,\chi_N^a,C_N^a): [0,T]^2 \to \reals^4$,
in the topology of uniform convergence on $[0,T]^2$. 
This implies the existence of limit points for 
this sequence. By Proposition \ref{prop-macro}
every such limit point is a solution of the integral
equations \req{eqC1}--\req{eqE} with the stated symmetry and 
boundary conditions. By Lemma \ref{lem-diff} each such solution 
results with $C$ and $\chi$ (i.e. $R$) 
that satisfy the integro-differential
equations \req{eqR}--\req{eqZ}. In view of Proposition \ref{uniqueness}
the latter system admits at most one solution per
given boundary conditions. Hence, we conclude that the sequence 
$(\chi_N^a,C_N^a)$ converges uniformly in $[0,T]^2$ to the unique solution
of \req{eqR}--\req{eqZ} subject to the appropriate boundary conditions.
Further, by \req{eq:asself} of Proposition \ref{self-average}
both $C_N-C_N^a$ and $\chi_N-\chi_N^a$ converge uniformly to zero,
almost surely. Thus, the solution of \req{eqR}--\req{eqZ}
is also the unique almost sure uniform (in $s,t$)
limit of $(\chi_N,C_N)$, as stated in Theorem \ref{thm-macro}. 
The $L_p$ convergence then follows by the uniform bounds on moments
of $C_N$ and $\chi_N$ (see Proposition \ref{tightchic}), 
thus completing the proof of the theorem.
\hfill\qed

\appendix\section{It\^o's calculus}\label{sec-ito}

Let $\left\lbc
\{x^i_t,y^i_t,z^i_t,w^i_t\}_{t\ge 0}, i\in\N\right\rbc$
be semi-martingales such that, 
$$d\langle r^i,p^j \rangle_t=\delta_{i=j} dt$$
for any $p,r\in\{x,y,z,w\}$.
Denoting, for $p,r\in\{x,y,z,w\}$, $s,t\ge 0$,
$N\in\N$,
$$K_{p,r}^N(s,t)
:= \oneN \sum_{i=1}^Np^i_sr^i_t\,,
$$
we already made use of the following simple
stochastic calculus lemma.
\begin{lem}\label{ito}
For any polynomials $P,Q$, and any $\s,\t,v\ge 0$,
\begin{eqnarray*}
P(K_{x,y}^N(\s,\t))Q(K_{z,w}^N(\s,v))
&=& P(K_{x,y}^N(0,\t))Q(K_{z,w}^N(0,v))\\
&&
+\int_0^{\s} P'(K_{x,y}^N(u,\t))Q(K_{z,w}^N(u,v)) d_uK_{x,y}^N(u,\t)\\
&&+ \int_0^{\s} P(K_{x,y}^N(u,\t))Q'(K_{z,w}^N(u,v)) d_uK_{z,w}^N(u,v)\\
&&+\frac{1}{2N}
K_{y,y}^N(\t,\t)\int_0^{\s}
P''(K_{x,y}^N(u,\t))Q(K_{z,w}^N(u,v)) du\\
&&+\frac{1}{2N}
K_{w,w}^N(v,v)\int_0^{\s}
P(K_{x,y}^N(u,\t))Q''(K_{z,w}^N(u,v)) du\\
&&+\oneN 
K_{y,w}^N(\t,v)\int_0^{\s} P'(K_{x,y}^N(u,\t))Q'(K_{z,w}^N(u,v)) du\\
\end{eqnarray*}
where
$$d_uK_{z,w}^N(u,v):=\oneN \sum_{i=1}^N w^i_v d z^i_u,$$
and all the stochastic integrals are defined via our convention
(of putting terms such as 
$y^{i_1}_{\t}\cdots y^{i_a}_{\t}  w^{j_1}_v\cdots w^{j_b}_v$ outside
the integral).
\end{lem}

\noindent
\prf By the bi-linearity of the formula given, it is enough to
prove the lemma for $P(x)=x^a$ and $Q(x)=x^b$. In this case, writing
$$(K_{x,y}^N(\s,\t))^a(K_{z,w}^N(\s,v))^b=
N^{-(a+b)} \sum_{i_1,\cdots i_a}\sum_{j_1,\cdots j_b } 
 y^{i_1}_{\t}\cdots y^{i_a}_{\t} w^{j_1}_v\cdots w^{j_b}_v 
\left(x^{i_1}_{\s}\cdots x^{i_a}_{\s} z^{j_1}_{\s}\cdots
z^{j_b}_{\s}\right)\, ,
$$
and using It\^o's formula for
$x^{i_1}_{\s}\cdots
x^{i_a}_{\s} z^{j_1}_{\s}\cdots
z^{j_b}_{\s}$ gives the stated result.
\qed

\section{Supremum of Gaussian processes indexed on large dimensional spheres.}
\label{sec-gauss}

In this section we prove the bound \req{eq:ger2} 
which is a direct consequence of the following general lemma about 
supremum of Gaussian processes indexed on large dimensional balls.
The outline for a direct proof of such a result 
by a chaining argument was kindly communicated to us by A. Bovier, 
whom we thank gratefully. This chaining argument can be adapted rather straightforwardly from \cite{BG}.
Anton Bovier also mentioned that this result should be the consequence of a more general one. 
We have indeed  found the proper way to see it as a consequence of classical and
well known facts on Gaussian processes, and to give simple references.
 
\nn

Let $(X_N(\bx))$ be a sequence of real valued,
centered Gaussian processes indexed by $\bx\in \R^{N}$.
Consider, for every $\rho>0$, the closed 
Euclidean ball $B_N(0,\rho)$ in  $\R^{N}$, and define 
$$X_N^*(\rho)= \sup_{\bx\in B_N(0,\rho)} \frac{|X_N(\bx)|} {\sqrt{N}}$$
We will also introduce the usual metric on $\R^{N}$ 
associated to the process $X_N$,
$$ d_X(\bx,\by)= \E[|X_N(\bx)-X_N(\by)|^2]^{1/2} $$
We denote by $\|\cdot\|$ the Euclidean norm and by $(\bx,\by)_N$ 
the corresponding inner product on $\R^{N}$.

\begin{lem}\label{Gaussian bound} 
Suppose that 
\begin{equation}\label{eq:ass1}
 \sup_{N} \E[ X_N(0)^2]  < \infty 
\end{equation}
and that
\begin {equation}\label{eq:ass2}
 \sup_{N}\sup_{\bx,\by\in B_N(0,\rho)} \frac{d_X(\bx,\by)}{\|\bx-\by\|} 
< \infty \,.
\end{equation}
Then, for every $k \in \N$
\begin{equation}\label{eq:res1}
\sup_{N} \E[X_N^*(\rho)^k] < \infty \,.
\end{equation}
Moreover, there exists a constant 
$\kappa<\infty$ such that for all $N$ and every $t>0$,
\begin{equation}\label{eq:res2}
\P[X_N^*(\rho) \geq \kappa +t] \le \exp(-N t^2/\kappa) \,.
\end{equation}
\end{lem}
\prf This result is a direct consequence of Dudley's theorem (\cite{Dudley}).
Indeed, the assumption \req{eq:ass2} implies that for any $N$ and
$\epsilon>0$ one can cover $B_N(0,\rho)$ by the union of  
certain $C(\rho) \epsilon^{-N}$ balls of radius $\epsilon$, 
in the metric $d_X$, where
the constant $C(\rho)$ depends on $\rho$ but not on the dimension $N$.
Thus, Dudley's theorem  (see also \cite[Theorem 11.17]{Ledoux-Talagrand})
shows that
$$ \E[ \sup_{\bx\in B_N(0,\rho)} X_N(\bx) ] \leq C'(\rho) \sqrt{N}\,, $$
where the constant $C'(\rho)$ again depends  only on $\rho$ and not on $N$.
Using the obvious fact that
\begin{eqnarray*}
\E[ \sup_{\bx\in B_N(0,\rho)} |X_N(\bx)| ] -
\E[|X_N(0)|]&\leq&
\E[ \sup_{\bx,\by\in B_N(0,\rho)}|X_N(\bx)-X_N(\by)|] \\
&=&
\E[ \sup_{\bx,\by\in B_N(0,\rho)}\{ X_N(\bx)-X_N(\by)\} ] 
= 2 \E[ \sup_{\bx\in B_N(0,\rho)} X_N(\bx) ]\,,
\end{eqnarray*}
and the assumption \req{eq:ass1}, we 
see that the conclusion \req{eq:res1} holds for $k=1$.
Thus, $X_N$ admits a version with almost all 
sample paths bounded and uniformly continuous on $B_N(0,\rho)$.
One can then consider $X_N$ as an (infinite dimensional) 
Gaussian vector in the space of 
continuous functions on the ball $B_N(0,\rho)$, equipped with 
the supremum norm.
It is also a well known fact that, 
for such a Gaussian vector, 
all moments of the norm are controlled by the first (e.g.
see the last statement of \cite[Corollary 3.2]{Ledoux-Talagrand}). This is 
thus enough to ensure that \req{eq:res1} holds
for every $k\in \N$.

The  tail estimate \req{eq:res2} is also
classical in the Gaussian context. 
For instance, the assumptions \req{eq:ass1} and \req{eq:ass2} immediately 
imply that the weak variance 
$$
\sigma(X_N)=\sup_{\bx \in B_N(0,\rho)} \Big\{ \E [ X_N(\bx)^2 ]^{1/2}
 \Big\}$$
of \cite[page 56]{Ledoux-Talagrand} is bounded in $N$.
 Hence,
by \cite[estimate (3.2), page 57]{Ledoux-Talagrand} it is easy to see that 
there exists a finite constant $\kappa > \sup_{N} \E[X_N^*(\rho)]$ 
for which \req{eq:res2} applies. 
\qed

We proceed to apply Lemma \ref{Gaussian bound} 
to the situation of interest here. To this end,
fixing $p\in\N$ consider the Gaussian process defined on $(\R^{N})^{p}$ by
$$X_{N,p}(\bx)=
\sum_{1\le i_j\le N, 1\le j\le p} G_{\{i_1,\cdots,i_p\}} x_{i_1}^1
x_{i_2}^2x_{i_3}^3\cdots x_{i_p}^p
$$
with independent, centered Gaussian variables $G_{\{i_1,\cdots,i_p\}}$
of variances $c(\{i_1,\ldots,i_p\})$ of \req{eq:vardef}.
Considering in $(\R^{N})^{p}$ the Cartesian product 
$$ 
B(\rho) = \prod_{i=1}^{p}{\{\bx^i \in \R^N, \|\bx^i\|\leq \rho}\}\,
$$
of $p$ Euclidean balls, 
we wish to estimate the moments and tail of
$\sup_{\bx\in B(\rho)} |X_{N,p}(\bx)|$.
To this end, note that since 
$X_{N,p}$ is a symmetric $p$-linear form on $\R^{N}$
(i.e., $X_{N,p}(\bx)$ is invariant to permutations of the 
vectors $\bx^1, \ldots, \bx^p$), we have by polarization that 
\begin{equation}\label{polar}
\sup_{\bx\in B(\rho)}|X_{N,p}(\bx)|\leq C(p) 
\sup_{\|{\bf u}\|\leq p\rho}\,|X_N({\bf u})|\,,
\end{equation}
where one defines $X_N({\bf u})$ for ${\bf u}\in \R^{N}$ by
$$ X_N(\buu)= X_{N,p}(\buu,\ldots,\buu)= 
\sum_{1\le i_j\le N} G_{\{i_1,\cdots,i_p\}} u_{i_1}
u_{i_2}u_{i_3}\cdots u_{i_p}\,.
$$

The Gaussian process $X_N(\buu)$ obviously satisfies
\req{eq:ass1} since $X_N(0)=0$ for all $N$. 
Turning to check that \req{eq:ass2} is satisfied as well,
note that
$$
X_N(\buu)=N^{(p-1)/2}\sum_{i=1}^N G^i(\buu) u_i \,,
$$
for the case where $\nu'(r)=(p-1)! r^{p-1}$. Thus, by
\req{eq:diagcov} and \req{eq:offdcov}, the
covariance of the process $X_N$ is 
$$\E[ X_N(\buu)X_N(\bv)] = p! (\buu,\bv)_N^p \,,$$
so that
$$ d_X(\buu,\bv)^2= p! \Big[
(\buu,\buu)_N^p +  (\bv,\bv)_N^p - 2(\buu,\bv)_N^p\Big] \,.
$$
It is then easy to check that for all $\buu$ and $\bv$,
$$ d_X(\buu,\bv) \leq C(p)\|\buu-\bv\|\,,$$
yielding that \req{eq:ass2} holds.
Hence, by Lemma \ref{Gaussian bound}, for every positive
integers $p$ and $k$,
\begin{equation}\label{bovier1}
\sup_{N}\E [\sup_{\bx\in B(\rho)}(\frac{|X_{N,p}(\bx)|}{\sqrt{N}})^k] 
\leq \sup_{N} \E[X_N^*(\rho)^k] < \infty \,.
\end{equation}
In view of \req{jnorm}, combining \req{bovier1} for $p=1,\ldots,m$ then
proves the bound \req{eq:ger2}.
We further have by \req{eq:res2} and \req{polar} that for some
$\wka < \infty$, any $p \leq m$, all $N$ and every $t>0$,
$$
\P[\sup_{\bx\in B(\rho)}(\frac{|X_{N,p}(\bx)|}{\sqrt{N}}) \geq 
\wka+t] \le \exp(-N t^2/\wka) \,,
$$
from which it follows that,
\begin{equation}\label{bovier2}
\P[\|\bJ\|_\infty^N \geq \wka+t] \le m \exp(-N t^2/\wka) \,.
\end{equation}

\section{Gaussian change of measure identities}\label{sec-alice}

The change of measure that is key to the proof of 
Lemma \ref{condexp} is a special case of the following ``standard''
Gaussian computation (compare \req{eq:lndf} and \req{eq:rn}).

\begin{prop}\label{prop-alice}
Suppose under the law $\P$ we have a
finite collection $\BJ = \{ J_\a \}_\a$ of 
non-degenerate, independent, centered Gaussian random variables,
and $G_s^i=\sum_\a J_\a L_s^i(\a)$ for 
all $s \in [0,\tau]$ and $i \leq N$, where for each $\a$
the coefficients $L_s^i$ which are independent of $\BJ$ 
are in $L_2(\{1,\ldots,N\}\ts [0,\tau])$. Suppose further that
$U_s^i$ is a continuous semi-martingale, independent of $\BJ$
and such that for each $\a$ the stochastic integral 
$$
\mu_\a := \sum_{i=1}^N \int_0^\tau L_u^i(\a) dU^i_u \,,
$$ 
is well defined and almost surely finite. Let 
$\P^*$ denote the law of $\BJ$ such that 
$\P^* = \L_\tau/\E(\L_\tau) \P$, where
\begin{equation}\label{eq:rn}
\L_\tau =\exp\Big\{ \sum_{i=1}^N \int_0^\tau G^i_s dU^i_s
-\half\sum_{i=1}^N \int_0^\tau (G^i_s)^2 ds\Big\}  \,.
\end{equation}
Let $k_{ts}^{ij}=\E(G_t^i G_s^j)$, $V_s^i=\E^*(G_s^i)$ and
$\Gamma_{ts}^{ij}=\E^*[(G_t^i-V_t^i)(G_s^j-V_s^j)]$. Then,
for any $s \leq \tau$ and $i \leq N$,
\begin{equation}\label{eq:mean}
V^i_s + [k_\tau V]^i_s = [k_\tau \circ dU]^i_s\,,
\end{equation}
and for any $s,t \leq \tau$ and $i,l \leq N$,
\begin{equation}\label{eq:cov}
\sum_{j=1}^N \int_0^\tau k_{su}^{ij} \Gamma_{ut}^{jl} du 
+ \Gamma_{st}^{il} = k_{st}^{il} \;.
\end{equation}
\end{prop}
\prf 
Let $v_\a=\E(J_\a^2)>0$ denote the variance of $J_\a$ and 
\begin{equation}\label{eq:chmes}
R_{\a\g} := \sum_{i=1}^N \int_0^\tau L_u^i(\a) L_u^i(\g) du \,,
\end{equation}
observing that  
$$
\L_\tau =\exp\Big\{ \sum_{\a} J_\a \mu_\a -
\half\sum_{\a,\g} J_\a J_\g R_{\a\g} \Big\} \,.
$$
With $\BD= {\rm diag} (v_\a)$ a positive definite matrix and
$\BR=\{R_{\a\g}\}$ positive semi-definite, it 
follows from this representation of $\L_\tau$
that under $\P^*$ the random vector $\BJ$ has a Gaussian law
with covariance matrix $(\BD^{-1}+\BR)^{-1}$ and mean vector 
$\bq=\{q_\a\}=(\BD^{-1}+\BR)^{-1} \bmu$. Hence, for any $\a$,
\begin{equation}\label{eqq}
q_\a + v_\a \sum_\g R_{\a\g} q_\g = v_\a \mu_\a \,.
\end{equation}
As $k_{su}^{ij}=\sum_\a L_s^i (\a) v_\a L_u^j (\a)$, it is not
hard to check that
$$
[k_\tau \circ dU]^i_s := \sum_{j=1}^N \int_0^\tau k_{su}^{ij} dU^j_u 
= \sum_\a L_s^i (\a) v_\a \mu_\a \;.
$$
Obviously, $V^i_s=\sum_\a L_s^i(\a) q_\a$ so we get \req{eq:mean}
out of \req{eqq} upon verifying that
$$
[k_\tau V]^i_s := \sum_{j=1}^N \int_0^\tau k_{su}^{ij} V^j_u du
= \sum_{\a,\g} L_s^i (\a) v_\a q_\g \sum_{j=1}^N 
\int_0^\tau L_u^j (\a) L_u^j(\g) du =
\sum_{\a,\g} L_s^i(\a) v_\a R_{\a\g} q_\g \,,
$$
with the last identity due to \req{eq:chmes}.
 
Turning to prove \req{eq:cov}, since $\Gamma_{ut}^{jl}$ is the 
covariance of $G_u^j$ and $G_t^l$ under the tilted law $\P^*$, 
we have that   
$$
\Gamma_{ut}^{jl} = \sum_{\a,\g} L_u^j (\a)[(\BD^{-1}+\BR)^{-1}]_{\a\g}
 L_t^l(\g) \,,
$$
and hence by \req{eq:chmes} we see that
$$
\sum_{j=1}^N \int_0^\tau k_{su}^{ij} \Gamma_{ut}^{jl} du 
= \sum_{\a,\g} L_s^i (\a) v_\a [\BR (\BD^{-1}+\BR)^{-1}]_{\a\g} L_t^l(\g) \;.
$$
With $\BD={\rm diag} (v_\a)$ we easily 
get \req{eq:cov} out of the 
matrix identity $({\bf I}+\BD \BR)(\BD^{-1}+\BR)^{-1} = \BD$.
\hfill\qed

\medskip 
{\bf Acknowledgement} 
We are gratefull to Leticia Cugliandolo and Jorge Korchan for 
their course at the 2001 Oberwolfach mini-workshop on aging,
which motivated this research. We also thank Anton Bovier 
for kindly communicating to us a direct proof of Lemma 
\ref{Gaussian bound}, and to SR Srinivasa Varadhan for 
helping us with the proof of Lemma \ref{loc-conc}.

\end{document}